\documentclass[11pt]{amsart}

\usepackage{enumerate}
\usepackage{latexsym}
\usepackage[pdftex]{graphicx}
\usepackage{amssymb}
\usepackage{xcolor}
\usepackage{stackrel}
\usepackage{systeme}
\usepackage[normalem]{ulem}

\newtheorem{theorem}{Theorem}
\newtheorem{lemma}{Lemma}
\newtheorem{proposition}{Proposition}
\newtheorem{remark}{Remark}
\newtheorem{example}{Example}

\newtheorem{definition}{Definition}
\newtheorem{corollary}{Corollary}

\newtheorem{problem}{Problem}
\newtheorem*{problem*}{Problem}

\def\R{{\mathbb R}}

\newcommand{\beq}{\begin{equation}}
\newcommand{\eeq}{\end{equation}}
\newcommand{\beqna}{\begin{eqnarray*}}
\newcommand{\eeqna}{\end{eqnarray*}}
\newcommand{\beqn}{\begin{equation*}}
\newcommand{\eeqn}{\end{equation*}}
\newcommand{\bp}{\begin{proof}}
\newcommand{\ep}{\end{proof}}
\newcommand{\bprop}{\begin{proposition}}
\newcommand{\eprop}{\end{proposition}}
\newcommand{\bt}{\begin{theorem}}
\newcommand{\et}{\end{theorem}}
\newcommand{\bex}{\begin{example}}
\newcommand{\eex}{\end{example}}
\newcommand{\bc}{\begin{corollary}}
\newcommand{\ec}{\end{corollary}}
\newcommand{\bl}{\begin{lemma}}
\newcommand{\el}{\end{lemma}}
\newcommand{\bprob}{\begin{problem}}
\newcommand{\eprob}{\end{problem}}
\newcommand{\br}{\begin{remark}}
\newcommand{\er}{\end{remark}}
\newcommand{\bd}{\begin{definition}}
\newcommand{\ed}{\end{definition}}

\begin{document}

\title
[On   bodies   floating in equilibrium in every orientation ]
{A negative answer to Ulam's Problem 19 from the Scottish Book}

\author[D. Ryabogin]{Dmitry Ryabogin}
\address{Department of Mathematical Sciences, Kent State University,
Kent, OH 44242, USA} \email{ryabogin@math.kent.edu}

\thanks{The   author is  supported in
part by   Simons Collaboration Grant for Mathematicians program 638576,
by U.S.~National Science Foundation Grant
DMS-1600753 and by United States - Israel Binational Science Foundation (BSF)}


\begin{abstract}
	We give a negative answer to   Ulam's Problem 19 from the Scottish Book  asking  {\it  is a solid of uniform density which will float in water in every position a sphere?}  Assuming that the density of water is $1$, we show  that  there exists a  strictly convex    body of revolution $K\subset {\mathbb R^3}$ of uniform density $\frac{1}{2}$, which is not a Euclidean ball,  yet 
 floats in equilibrium in every orientation. We prove an analogous result in all dimensions $d\ge 3$.
\end{abstract}

\maketitle

\section{Introduction}

The following intriguing problem was  proposed by Ulam \cite[Problem 19]{U}: {\it If a convex body  $K\subset {\mathbb R^3}$ made of material of uniform density ${\mathcal D}\in(0,1)$  floats in equilibrium in any orientation $\textnormal{(}$in water, of density $1$$\textnormal{)}$, must $K$ be a Euclidean ball? }

Schneider \cite{Sch1} and
 Falconer   \cite{Fa}  showed that this is true, provided $K$ is centrally symmetric and ${\mathcal D}=\frac{1}{2}$.  
No results are known for other densities ${\mathcal D}\in (0,1)$ and no counterexamples have been found so far.

The ``two-dimensional version" of the problem is also very interesting. In this case, we consider  floating logs of uniform cross-section, and seek for the ones that will float in every orientation with the axis horizontal.
If ${\mathcal D}=\frac{1}{2}$, Auerbach  \cite{A} has exhibited logs with non-circular cross-section, both convex and non-convex, whose  boundaries  are   so-called Zindler curves \cite{Zi}.
More recently,    Bracho,  Montejano and   Oliveros  \cite{BMO} showed that 
 for densities ${\mathcal D}=\frac{1}{3}$, $\frac{1}{4}$, $\frac{1}{5}$ and  $\frac{2}{5}$ the answer is affirmative, while  Wegner 
 proved that for some  other  values of ${\mathcal D}\neq \frac{1}{2}$ the answer is negative, \cite{Weg1}, \cite{Weg2}; see also related  results of V\'arkonyi \cite{V1}, \cite{V2}. Overall, the case of general ${\mathcal D}\in (0,1)$ is notably involved and widely open.

  In this paper  we prove the following result.
\bt\label{bitsya}
Let $d\ge 3$. There exists   a  strictly convex  \textnormal{non-centrally-symmetric}  body of revolution $K\subset {\mathbb R^d}$  which  floats in equilibrium in every orientation at the level $\frac{\textnormal{vol}_d(K)}{2}$.
\et
This gives 
 \bt\label{hob3}
 The answer to Ulam's Problem 19 is negative, i.e., there exists a convex body $K\subset {\mathbb R^3}$   of density ${\mathcal D}=\frac{1}{2}$,  which is not a Euclidean ball, yet  floats in equilibrium in every orientation.
 \et

Our bodies will be {\it small perturbations of the Euclidean ball}.   We combine our recent results from \cite{R} together with  work of Olovjanischnikoff  \cite{O},
and then use the machinery developed together with Nazarov and Zvavitch in \cite{NRZ}. The proofs of Theorem \ref{bitsya} for even and odd $d$ are different. 
For even $d$  we solve a finite moment problem to obtain our body as a {\it local} perturbation of the Euclidean ball.
The case  $d\ge 3$ with odd $d$ is  more involved. To control the perturbation, we use the properties of the  spherical Radon transform,  \cite[pp. 427-436]{Ga}, \cite[Chapter III, pp. 93-99]{He}.

We refer the reader to   \cite[pp. 19-20]{CFG},  \cite[pp. 376-377]{Ga}, \cite{G},  \cite[pp. 90-93]{M} and \cite{U}  for an exposition of known results related to the problem.

This paper is structured as follows. In Section \ref{next1}, we recall all the necessary notions and statements needed to  prove the  main result.  In Section \ref{BRev}, we reduce the problem to finding  a non-trivial solution to a  system of two integral equations.  In Section \ref{Even}, we prove Theorem \ref{bitsya} for even $d$. In Section \ref{Odd}, we give the proof of  
Theorem \ref{bitsya} for odd  $d$ and prove Theorem \ref{hob3}.
 In Appendix A, we present the proof of Theorem \ref{olovzhal} given in \cite{O}. We prove the converse part of Theorem \ref{Fedja1}  in Appendix B.

\section{Notation and auxiliary results}\label{next1}

Let ${\mathbb N}=\{1,2,\dots,\}$ be the set of natural numbers.
A convex body $K\subset {\mathbb R^d}$, $d\ge 2$,   is a convex compact set with  non-empty interior $\textrm{int} K$. The boundary of $K$ is denoted by
 $\partial K$. We say that $K$ is  strictly convex  
  if $\partial K$ does not contain a segment.
 We say that $K$ is origin-symmetric if $K=-K$ and centrally-symmetric if there exists $p\in{\mathbb R^d}$ such that $K-p=\{q-p:\,q\in K\}$ is origin-symmetric.
 Let $S^{d-1}=\{\xi\in{\mathbb R^d}:\,\sum\limits_{j=1}^d\xi_j^2= 1\}$  be the unit sphere in ${\mathbb R^d}$ centered at the origin and  let $B_2^d=\{p\in{\mathbb R^d}:\,\sum\limits_{j=1}^dp_j^2\le 1\}$ be the unit Euclidean ball centered at the origin.
 We denote by $\kappa_d=\textrm{vol}_d(B^d_2)$
 the $d$-dimensional volume of $B^d_2$ and we let $e_1,\dots, e_d$ be
 the standard basis in ${\mathbb R^d}$.
Given  $\xi\in S^{d-1}$, we denote by
$\xi^{\perp}=\{p\in{\mathbb R^d}:\, p\cdot \xi=0\}$ the    subspace orthogonal to $\xi$, where   $p\cdot\xi=p_1\xi_1+\dots +p_d\xi_d$ is the usual inner product in ${\mathbb R^d}$.  
For $p\in {\mathbb R^d}$ we put $|p|=\sqrt{p_1^2+\dots+p_d^2}$.
We also denote by 
${\mathcal B}(\xi, \rho)=\{ p\in S^{d-1}:\,p\cdot \xi>\rho\}$ 
the spherical cap centered at $\xi\in S^{d-1}$ of radius $\rho\in [-1,1)$; we tacitly assume that ${\mathcal B}(\xi, -1)=S^{d-1}$.
We say that a  hyperplane $H$ is the supporting hyperplane of  a convex body $K$ if $K\cap H\neq \emptyset$, but $\textrm{int}K\cap H=\emptyset$. Let $W_j$ be a $j$-dimensional  plane in ${\mathbb R^d}$, $1\le j\le d$.
The {\it center of mass} of a  compact convex set $K\subset W_j$ with a non-empty relative interior will be denoted by 
${\mathcal C}(K)=\frac{1}{\textrm{vol}_j(K)}
\int\limits_{K}pdp$,
where $\textrm{vol}_j(K)$ is the $j$-dimensional volume of $K$ and $dp$ stands for the usual Lebesgue measure in ${\mathbb R^j}$.  Given two sets $A$ and $B$ in ${\mathbb R^d}$, we denote by $A\times B$ their Cartesian product, i.e., the set of ordered pairs $\{(a,b):\,a\in A, b\in B\}$.
Let  $k\in{\mathbb N}$.  We say that a function
$h: {\mathbb R}\to{\mathbb R}$ supported on  a closed interval $[a,b]\subset {\mathbb R}$, $a<b$, is in $C^k$ (in $C^{\infty}$) if it has continuous derivatives up to  order $k$ (of all orders). We define  its norm as
$\|h\|_{C^k}=\sum\limits_{m=0}^ {k}\max\limits_{\{s\in [a,b]\}}|h^{(m)}(s)|$,
where $h^{(m)}$ is the $m$-th derivative of $h$.
We say that  a convex body $K\subset{\mathbb R^d}$ is of class $C^k$  if $K$ has a $C^k$-smooth  boundary, i.e.,   for every point $z\in \partial K$  there exists a neighborhood $U_z$ of $z$  in ${\mathbb R^d}$ such that  $\partial K\cap U_z$ can be written as a graph of a function having 
all continuous partial derivatives  up to  the $k$-th order.

\begin{figure}[h]
	\centering
	\includegraphics[height=3.3in]{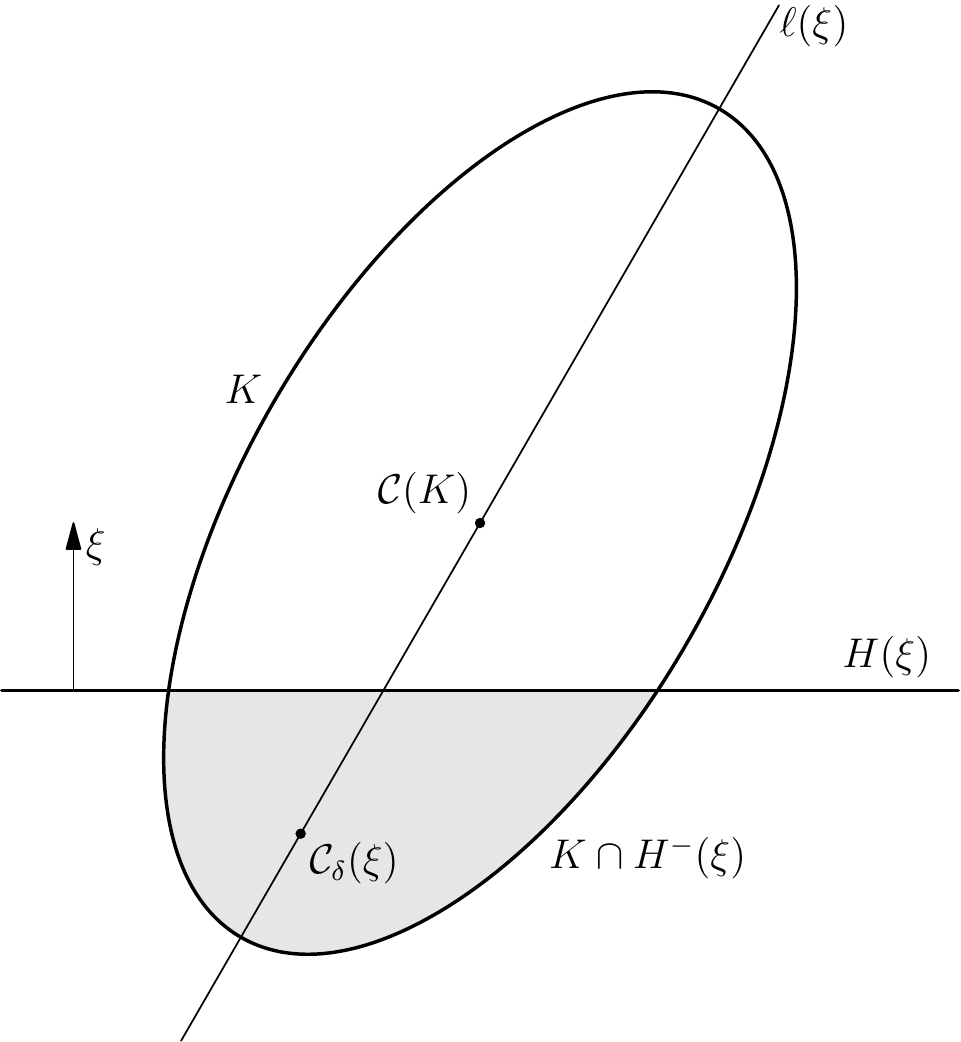} 
	\caption{A body $K$, its submerged part  $K\cap H^-(\xi)$ and the line $\ell(\xi)$ passing through ${\mathcal C}(K)$ and ${\mathcal C}_{\delta}(\xi)$.}
	\label{width1}
\end{figure}

Let $d\ge 3$, let $K\subset {\mathbb R^d}$ be a convex body 
and let $\delta\in (0, \textnormal{vol}_d(K))$ be fixed. Given a direction $\xi\in S^{d-1}$ and $t=t(\xi)\in{\mathbb R}$, we call a hyperplane 
$$
H(\xi)=\{p\in{\mathbb R^d}:\,p\cdot\xi=t\}
$$
the {\it cutting hyperplane} of $K$ in the direction $\xi$,  if it cuts out of $K$ the given volume $\delta$, 
i.e., if 
\begin{equation}\label{fubu}
\textrm{vol}_d(K\cap H^-(\xi))=\delta,\qquad\textnormal{where}\quad H^-(\xi)=\{p\in{\mathbb R^{d}}:\,p\cdot\xi\le t(\xi)\},
\end{equation}
(see Figure \ref{width1}).

We recall several well-known facts and definitions (see \cite[Chapter XXIV]{DVP}, \cite[Chapter 2]{L}, \cite[Chapter 4]{Tu}, \cite[Hydrostatics, Part I]{Zh}). 
\bd\label{efb}
Let   $\xi\in S^{d-1}$  and let ${\mathcal C}_{\delta}(\xi)$ be the center of mass of the submerged part $K\cap H^-(\xi)$ satisfying (\ref{fubu}). 
We say that $K$ {\it floats in equilibrium in the direction $\xi$ at the level $\delta$} if  the line $\ell(\xi)$ passing through ${\mathcal C}(K)$ and ${\mathcal C}_{\delta}(\xi)$ is orthogonal to the ``free water surface" $H(\xi)$, i.e., the line $\ell(\xi)$ is ``vertical" $\textnormal{(}$parallel to $\xi$, see Figure \ref{width1}$\textnormal{)}$.
We say that $K$ floats in equilibrium in every orientation at the level $\delta$ if $\ell(\xi)$ is parallel to $\xi$ for every $\xi\in S^{d-1}$.
\ed

\bd\label{scn}
The geometric locus  $\{{\mathcal C}_{\delta}(\xi):\,\xi\in S^{d-1}\}$  is called the {\it surface of centers} ${\mathcal S}={\mathcal S}_{\delta}$ or the {\it  surface of buoyancy}. 
\ed

Now  we recall the notion of  {\it  characteristic points} of  a family of hyperplanes (cf. 
\cite[pp. 107-110]{BG}, \cite[pp. 48-50]{Wea}, or \cite[pp. 26-54]{Za}).  
\bd\label{chp1}
Let $d\ge 2$, $\xi_0\in S^{d-1}$, and let $\rho\in [-1,1)$. Consider    a family ${\mathcal Q}$ of hyperplanes in ${\mathbb R^d}$ 
such that for every direction $\xi\in {\mathcal B}(\xi_0, \rho)$ there exists a   hyperplane in $ {\mathcal Q}$ orthogonal to $\xi$. Assume also that for any $H\in {\mathcal Q}$, for any $(d-2)$-dimensional subspace $\Gamma$ parallel to $H$ and for any sequence $\{H_k\}_{k=1}^{\infty}$  of hyperplanes $H_k\in {\mathcal Q}$ converging to $H$ as $k\to\infty$ and parallel to $\Gamma$, the limit 
$\lim\limits_{k\to\infty}H\cap H_k$ exists.
 Given $H\in {\mathcal Q}$, we call a point $e\in H$  the  characteristic point of $ {\mathcal Q}$ with respect to $H$
if for any  $\Gamma$ and $\{H_k\}_{k=1}^{\infty}$, as above, we have $e\in \lim\limits_{k\to\infty}H\cap H_k$.
\ed

We will  need the following result from \cite{O}  (see  the lemma on pp. 114-117 and Remark 1 on p. 117).
\bt\label{olovzhal}
Let $d\ge 3$, let $K\subset {\mathbb R^d}$ be a convex body and let $\delta\in(0, \textnormal{vol}_d(K))$.
The characteristic points of the family   of  cutting hyperplanes $\{H(\xi):\xi\in S^{d-1}\}$  for which (\ref{fubu}) holds are the centers of mass of the  sections $\{K\cap H(\xi):\, \xi\in S^{d-1}\}$. 

Conversely, if the characteristic points of the  family of hyperplanes $\{H(\xi):\,\xi\in S^{d-1}\}$  intersecting the interior of $K$ and corresponding to the sections  $\{K\cap H(\xi):\,\xi\in S^{d-1}\}$ coincide with the centers of mass of
these sections, then the function $\xi\mapsto \textnormal{vol}_d(K\cap H^-(\xi))$ is  constant  on $S^{d-1}$ and   the constant is equal  to some $\delta \in(0, \textnormal{vol}_d(K))$.
\et
Since the reference \cite{O} is not readily available, for  convenience of the reader we present  the proof  of Theorem \ref{olovzhal}  in Appendix A.

To define the {\it moments of inertia} (see \cite[p. 553]{Zh}),
  consider a convex body $K$ and a
hyperplane $H(\xi)$ for which    (\ref{fubu}) holds.
Choose  any $(d-2)$-dimensional plane $\Pi\subset H(\xi)$ passing through the center of mass ${\mathcal C}(K\cap H(\xi))$ and let 
$\eta_1,\dots,\eta_{d-2}, \eta_{d-1}$ be an orthonormal basis of $\xi^{\perp}=\{p\in{\mathbb R^d}:\,p\cdot\xi=0  \}$ such that 
\begin{equation}\label{base}
\Pi={\mathcal C}(K\cap H(\xi))+\textrm{span}(\eta_1,\dots,\eta_{d-2}),\quad H(\xi)={\mathcal C}(K\cap H(\xi))+\xi^{\perp}. 
\end{equation}

\begin{figure}[h]
	\centering
	\includegraphics[height=2in]{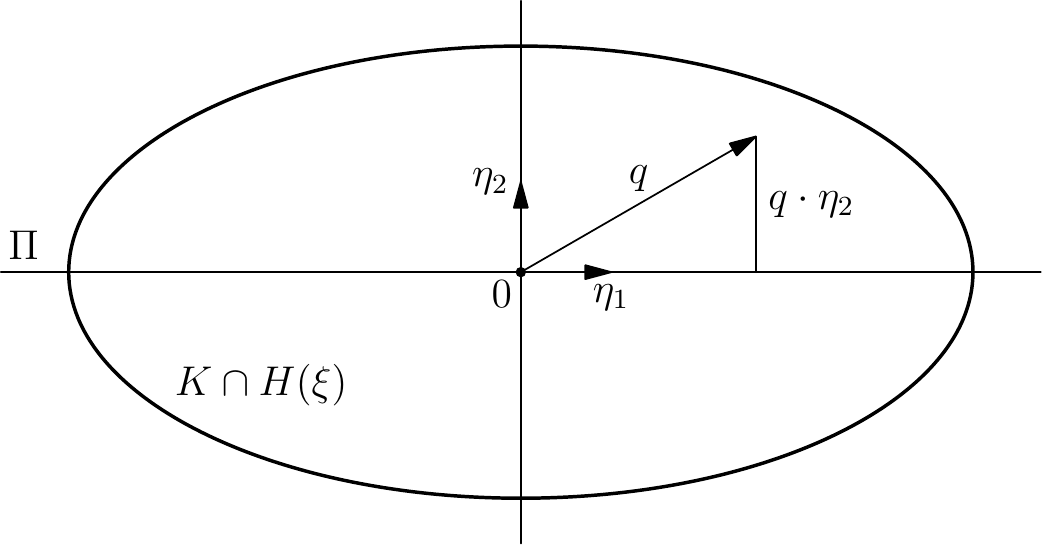} 
	\caption{Two-dimensional body $K\cap H(\xi)$ with  center of mass  at the origin, and a line $\Pi$ parallel to $\eta_1$; we have $\textnormal{dist}(\Pi,q)^2=|q|^2-(q\cdot \eta_1)^2=(q\cdot \eta_2)^2$.}
	\label{width4}
\end{figure}

\bd\label{miab}
The moment of inertia $I_{K\cap H(\xi)}(\Pi)$ of $K\cap H(\xi)$ with respect to  $\Pi$  is calculated by summing $\textnormal{dist}(\Pi,p)^2$ for every ``particle" $p$ in the set $K\cap H(\xi)$, 
$\textnormal{(}$see Figure \ref{width4}$\textnormal{)}$,
i.e., 
\begin{equation}\label{moment}
I_{K\cap H(\xi)}(\Pi)=\int\limits_{K\cap H(\xi)}\textnormal{dist}(\Pi,p)^2dp=
\int\limits_{K\cap H(\xi)-{\mathcal C}(K\cap H(\xi))}(q\cdot\eta_{d-1})^2\,dq,
\end{equation}
where $\textnormal{dist}(\Pi,p)=\min\limits_{\{q\in \Pi\} }|p-q|$.
\ed

We will  use the converse part of the following  theorem (see \cite[Theorem 1]{R} or \cite[Theorem 1.1]{FSWZ}\footnote{
  It is assumed in \cite{FSWZ} that 
in the case $\delta= \frac{\textrm{vol}_d(K)}{2}$
 the set of characteristic points of the cutting hyperplanes  is  a single point.}).

\bt\label{Fedja1}
Let  $d\ge 3$, let $K\subset{\mathbb R^d}$ be a convex body  and let $\delta\in (0, \textnormal{vol}_d(K))$.

If $K$ floats in equilibrium at the level $\delta$ in every orientation, then for all $\xi\in S^{d-1}$ 
and for all 
$(d-2)$-dimensional planes  $\Pi\subset H(\xi)$ passing through the center of mass ${\mathcal C}(K\cap H(\xi))$,
the cutting sections $K\cap H(\xi)$ have equal  moments of inertia  independent of $\xi$ and $\Pi$.

	Conversely, let   ${\mathcal C}({\mathcal S})={\mathcal C}(K)$.
	If for all cutting hyperplanes $H(\xi)$, $\xi\in S^{d-1}$, and for all 
	$(d-2)$-dimensional planes $\Pi\subset H(\xi)$ passing through the center of mass ${\mathcal C}(K\cap H(\xi))$,
	the cutting sections $K\cap H(\xi)$  have equal  moments of inertia independent of $\xi$ and $\Pi$, then  $K$ floats in equilibrium in every orientation at the level $\delta$.
\et

 For convenience of the reader we prove  the converse part of this theorem  in Appendix B\footnote{It is assumed in  \cite{R}  that $K$ is of class $C^1$. We give a slightly different proof that does not use this assumption.    }. 

\br\label{rv}
Let $\delta=\frac{\textnormal{vol}_d(K)}{2}$. Since for any 
$\xi\in S^{d-1}$, ${\mathcal C}(K)$ is the arithmetic average of ${\mathcal C}(K\cap H^+(\xi))$ and ${\mathcal C}(K\cap H^-(\xi))$,
 the condition ${\mathcal C}({\mathcal S})={\mathcal C}(K)$ is  satisfied and ${\mathcal S}$ is symmetric with respect to ${\mathcal C}(K)$.
\er

\section{Reduction to a system of integral equations}\label{BRev}

\begin{figure}[ht]
	\includegraphics[width=300pt]{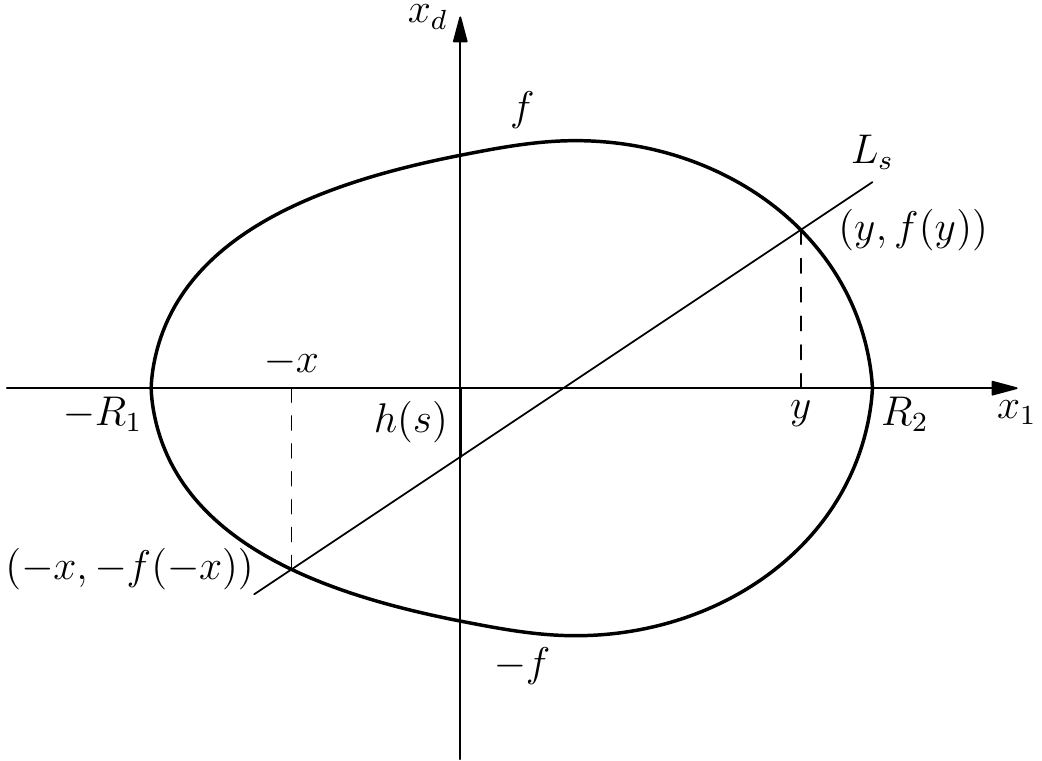}
	\caption{Sections of $K_f$ and $H(L)$  by the  $(x_1,x_d)$-plane.}
	\label{pic0}
\end{figure}

Let $d\ge 3$. We follow the notation from  \cite{NRZ}.
We will be dealing with bodies of revolution
$$
K_f=\{x\in{\mathbb R^d}:\,x_2^2+x_3^2+\dots+x_d^2\le f^2(x_1)  \}
$$
obtained by the rotation of a smooth concave function supported on $[-R_1,R_2]$ about the $x_1$-axis. 
Let $L(s,t)=L_s(t)=st+h(s)$ be a linear function with slope $s\in {\mathbb R}$, and let 
$$
H(L_s)=\{x\in{\mathbb R^d}:\,x_d=L_s(x_1)\}
$$ 
be the corresponding hyperplane.  The function $h$ will be chosen later. For now it is enough to assume that it   is  infinitely smooth,  not  identically  zero,  supported  on  $[1-2\tau,1-\tau]$   for  some  small $\tau>0$, and $h$ and sufficiently many of its derivatives are small. 
Let $-x=-x(s)$ and $y=y(s)$ be the first coordinates of the points of intersection of $\pm f$ and $L$
(see Figure \ref{pic0}).

To construct a system of two integral equations  we will prove four lemmas.
Consider the family of hyperplanes
\begin{equation}\label{fam1}
{\mathcal F}=\{H(L_s):\, s\in [0,\infty)\}.
\end{equation}

\bl\label{charp2}
Let $E$ be the set of  characteristic points of 
${\mathcal F}$. Then,
\begin{equation}\label{char11}
E=\{(-h'(s),0,\dots,0, L(s,-h'(s)))\in{\mathbb R^d}:\,s\in [0,\infty)\}.
\end{equation}
\el
\bp
Let ${\mathcal G}$ be the family of lines
${\mathcal G}=\{\ell_s:\,s\in[0,\infty)\}$, where each line $\ell_s$  is the intersection of $H(L_s)$ and  the $x_1x_d$-plane.
It  is enough  to show that 
$$
E\cap \{x_1x_d\!-\!\textnormal{plane}\}=\{(-h'(s), -sh'(s)+h(s))\in{\mathbb R^2}:\,s\in [0,\infty)\}.
$$
 We will use Definition \ref{chp1}. 
Let $s\in 
(1-2\tau, 1-\tau)$ and let $\ell_s\in {\mathcal G}$. Choose any sequence $\{\ell_{s_k}\}_{k=1}^{\infty}$, $\ell_{s_k}\in {\mathcal G}$,   converging
to $\ell_{s}$ as $k\to\infty$, and let $\{u_{s_k}\}_{k=1}^{\infty}$, $\{u_{s_k}\}=\ell_{s}\cap\ell_{s_k}$, be the corresponding sequence of points of intersection. Solving the system of two linear equations we see that $$
u_{s_k}=\Big(\frac{h(s)-h(s_k)}{s_k-s}, s_k\frac{h(s)-h(s_k)}{s_k-s}+h(s_k)\Big). 
$$
Hence, $\lim\limits_{s_k\to s}u_{s_k}$ exists and the point
$\lim\limits_{s_k\to s}u_{s_k}=(-h'(s), -sh(s)+h(s))$
is the characteristic point of ${\mathcal G}$ with respect to $\ell_s$. 

Next, we observe that $(0,0)$ is the characteristic point of ${\mathcal G}$ with respect to $\ell_{1-2\tau}$. Indeed,  it is enough to choose two sequences  of lines in ${\mathcal G}$, $\{    \ell_{s_k}   \}_{k=1}^{\infty}$, $\{\ell_{s_k'}\}_{k=1}^{\infty}$, both converging to $\ell_{1-2\tau}$, such that $s_k\in (1-2\tau, 1-\tau)$ and $s_k'\in (0,1-2\tau)$, and to use the fact  that  $\ell_{s_k'}\cap \ell_{1-2\tau}=
\{(0,0)\}$ for any line $\ell_{s_k'}$ with $s_k'\in (0,1-2\tau)$. Similarly, to show that
$(0,0)$ is the characteristic point of ${\mathcal G}$ with respect to $\ell_{1-\tau}$,
it is enough to choose the corresponding sequences $\{\ell_{s_k}\}_{k=1}^{\infty}$, $\{\ell_{s_k''}\}_{k=1}^{\infty}$, both converging to $\ell_{1-\tau}$, where $s_k\in (1-2\tau, 1-\tau)$ and $s_k''\in (1-\tau,\infty)$.

To finish the proof, it remains to observe that since  $h$ is supported by $[1-2\tau, 1-\tau]$, any two lines $\ell_s$, $\ell_{s'}$,
$s,s'\in 
[0,1-2\tau)\cup(1-\tau,\infty)$,
 intersect at $(0,0)$. Hence, $(0,0)$ is the characteristic point of ${\mathcal G}$ with respect to  any line $\ell_s$ for
$s\in 
[0,1-2\tau]\cup[1-\tau,\infty)$.
\ep
\bl\label{FL}
Let $s>0$. The condition 
\begin{equation}\label{fam11}
{\mathcal C}(K_f\cap H(L_s))=(-h'(s),0,\dots,0, L(s,-h'(s)))
\end{equation}
 reads as 
\begin{equation}\label{cmass}
\int\limits_{-x(s)}^{y(s)}(t+h'(s))(f(t)^2-L(s,t)^2)^{\frac{d-2}{2}}dt=0.
\end{equation}

Let 
$$
\Pi_1=\{x\in H(L_s):\,x_1=-h'(s)\},\qquad\quad \Pi_j=\{x\in H(L_s):\,x_j=0\},
$$ 
 $
j=2,\dots, d-1$.
The  moments of inertia  conditions 
$$
I_j=I_{K_f\cap H(L_s)}(\Pi_j)=const,\quad j=1,\dots, d-1,
$$
  read as
\begin{equation}\label{i1}
I_1=\kappa_{d-2}(1+s^2)^{\frac{3}{2}}\int\limits_{-x(s)}^{y(s)}(t+h'(s))^2(f(t)^2-L(s,t)^2)^{\frac{d-2}{2}}dt=const,
\end{equation}
\begin{equation}\label{ij}
I_j=\gamma_{d-2}\,\sqrt{1+s^2}\,\int\limits_{-x(s)}^{y(s)}(f(t)^2-L(s,t)^2)^{\frac{d}{2}}dt=const,
\end{equation}
where
$$
\gamma_{d-2}=\int\limits_{B^{d-2}_2}p_j^2dp,\quad j=2,\dots,d-1.
$$
\el
\bp
Fix $s>0$. Observe that the slice $K_f\cap H(L_s)\cap H_{t}$ of the cutting section $K_f\cap H(L_s)$ by the hyperplane $H_{t}=\{x\in{\mathbb R^d}:\, x_1=t\}$, $-x(s)<t<y(s)$, is the $(d-2)$-dimensional Euclidean ball
$$
B_2^{d-2}((t,0,\dots, 0,L(s,t)),r)\!=\!\{(t, x_2,\dots, x_{d-1}, L(s,t)): x_2^2+\dots+x_{d-1}^2\le r^2\}
$$
of radius $r=\sqrt{f^2(t)-L^2(s,t)}$ centered at $(t,0,\dots, 0,L(s,t))$. Hence,  for the  first coordinate of the center of mass in (\ref{fam11}) we have
\begin{equation}\label{mults}
\int\limits_{-x(s)}^{y(s)}(t+h'(s))
dt\int\limits_{B_2^{d-2}((t,0,\dots, 0,L(s,t)),r)}dp=
\end{equation}
$$
\kappa_{d-2}\int\limits_{-x(s)}^{y(s)}(t+h'(s))
(f(t)^2-L(s,t)^2)^{\frac{d-2}{2}}dt=0.
$$
This gives (\ref{cmass}). 

Similarly, since the distance in $K_f\cap H(L_s)$ between the points  $(t,x_2,\dots, x_d)$ $\in K_f\cap H(L_s)\cap H_{t}$ 
and
$(-h'(s),x_2,\dots, x_d)\in K_f\cap H(L_s)\cap H_{-h'(s)}$  is $\sqrt{1+s^2}|t+h'(s)|$, we have
$$
I_1=\sqrt{1+s^2}\int\limits_{-x(s)}^{y(s)}(\sqrt{1+s^2}(t+h'(s))^2dt
\int\limits_{B_2^{d-2}((t,0,\dots, 0,L(s,t)),r)}dp=
$$
$$
\kappa_{d-2}(1+s^2)^{\frac{3}{2}}\int\limits_{-x(s)}^{y(s)}(t+h'(s))^2(f(t)^2-L(s,t)^2)^{\frac{d-2}{2}}dt,
$$
proving (\ref{i1}).
Finally,  the expression in the left-hand side of  (\ref{ij}) for the other  moments can be obtained as
$$
I_j=\sqrt{1+s^2}\int\limits_{-x(s)}^{y(s)}dt
\int\limits_{B_2^{d-2}((t,0,\dots, 0,L(s,t)),r)}p_j^2dp=
$$
$$
\sqrt{1+s^2}\,\gamma_{d-2}\,\int\limits_{-x(s)}^{y(s)}(f(t)^2-L(s,t)^2)^{\frac{d}{2}}dt.
$$
\ep
\bl\label{work1h}
Let $s_o\ge 0$, let $K_f$ be as above and let ${\mathcal F}$ be the family of hyperplanes
defined as in (\ref{fam1}) for $s\ge s_0$, so that (\ref{fam11}) holds for $s\ge s_o$. Then for all $s>s_o$ and 
for all 
$(d-2)$-dimensional planes $\Pi\subset H(L_s)$ passing through the center of mass ${\mathcal C}(K_f\cap H(L_s))$,
the cutting sections $K_f\cap H(L_s)$  have equal   moments of inertia $I_{K_f\cap H(L_s)}(\Pi)$ independent of $s$ and $\Pi$,
provided (\ref{i1}) and (\ref{ij}) hold with  the same  constant  on the right-hand side, which is  independent of $s$ and $j=1,\dots,d-1$.
\el
\bp
Let $s_o\ge 0$ and    let $s>s_o$ be fixed. If $\Pi\subset H(L_s)$
is any $(d-2)$-dimensional plane  passing through the center of mass ${\mathcal C}_s={\mathcal C}(K_f\cap H(L_s))$, then by (\ref{moment}) we have
$$
I_{K_f\cap H(L_s)}(\Pi)=
\int\limits_{K_f\cap H(L_s)}((u-{\mathcal C}_s)\cdot\eta)^2\,du,
$$
where $\eta=\eta_{d-1}$ is a unit vector in the hyperplane  $H(L_s)-{\mathcal C}_s$ which is orthogonal to $l$.  

Let $\iota_1,\dots\iota_{d-1}$ be the orthonormal basis in $H(L_s)-{\mathcal C}_s$
such that $\iota_1\in\textrm{span}\{e_1,e_d\}$ and 
$\iota_j=e_j$ for $j=2,\dots, d-1$.
Decomposing  $\eta$ in this basis as $\sum\limits_{j=1}^{d-1}{\eta}_{(j)}\iota_j$,
we have
$$
I_{K_f\cap H(L_s)}(\Pi)=\sum\limits_{j=1}^{d-1}\eta_{(j)}^2\int\limits_{K_f\cap H(L_s)}((u-{\mathcal C}_s)\cdot \iota_j)^2\,du\,+
$$
$$
\sum\limits_{{\substack{j,l=1\\ j\neq l}}}^{d-1}
\eta_{(j)}\eta_{(l)}\int\limits_{K_f\cap H(L_s)}((u-{\mathcal C}_s)\cdot \iota_j)((u-{\mathcal C}_s)\cdot \iota_l)\,du=J_1+J_2.
$$
Using  the fact that $\eta$ is a unit vector, together with   (\ref{i1}) and (\ref{ij}), we have that $J_1$ is constant.

  We claim that $J_2=0$. Indeed, if $j$  is equal to $1$, then arguing as in the previous lemma, and 
  using the fact that
  $\int\limits_{B^{d-2}_2}p_ldp=0$ for  $l=2,\dots,d-1$,
  we see that
  $$
  \int\limits_{K_f\cap H(L_s)}((u-{\mathcal C}_s)\cdot \iota_1)((u-{\mathcal C}_s)\cdot \iota_l)\,du=
  $$
  $$
  \sqrt{1+s^2}\int\limits_{-x(s)}^{y(s)}(t+h'(s))
  dt\int\limits_{B_2^{d-2}((t,0,\dots, 0,L(s,t)),r)}p_ldp=0.
  $$
  The case when $l=1$ is similar.
  
If $j\neq 1$, $l\neq 1$, then we use the fact that
$\int\limits_{B^{d-2}_2}p_jp_ldp=0$ for  $j,l=2,\dots,d-1$, $j\neq l$, to obtain 
$$
\int\limits_{K_f\cap H(L_s)}((u-{\mathcal C}_s)\cdot \iota_j)((u-{\mathcal C}_s)\cdot \iota_l)\,du=
\int\limits_{-x(s)}^{y(s)}
dt\int\limits_{B_2^{d-2}((t,0,\dots, 0,L(s,t)),r)}p_jp_ldp=0.
$$
Thus, $I_{K_f\cap H(L_s)}(\Pi)$ is a constant  independent of $s$ and  of the  arbitrarily chosen $\Pi$.
The lemma is proved.
\ep
\bl\label{cl1}
Let $s_o\ge 0$. Assume that  (\ref{cmass}) is valid 
 for all $s>s_o$. 
 Then  (\ref{ij}) holds for all $s>s_o$ with
  the constant  independent of $s$ if and only if  (\ref{i1}) holds  for all $s>s_o$ with the constant independent of $s$. 
\el
\bp
  We recall that 
 \begin{equation}\label{defxy}
 \,\,\,L(s,t)=st+h(s),\quad f(y(s))=L(s, y(s)),\quad f(-x(s))=
 L(s, -x(s)),
 \end{equation}
 for $s\in {\mathbb R}$. Let $s_o\ge 0$ and let $s>s_o$. 
 We rewrite (\ref{ij}) as 
 $$
\int\limits_{-x(s)}^{y(s)}(f(t)^2-L(s,t)^2)^{\frac{d}{2}}dt=\frac{const}{ \gamma_{d-2}\, \sqrt{1+s^2}}
 $$
and  differentiate both sides 
 with respect to $s$ using (\ref{defxy}).  We have
 $$
 \int\limits_{-x(s)}^{y(s)}(f(t)^2-L(s,t)^2)^{\frac{d-2}{2}}(st+h(s))(t+h'(s))dt=\frac{const\,\,s}{ d\gamma_{d-2}\, (1+s^2)^{\frac{3}{2}}}.
 $$
Adding and subtracting $sh'(s)$ in the second parentheses under the integral  and using (\ref{cmass}), the last equality yields
$$
 s\,\int\limits_{-x(s)}^{y(s)}(f(t)^2-L(s,t)^2)^{\frac{d-2}{2}}(t+h'(s))^2dt=\frac{const\,\,s}{ d\gamma_{d-2}\, (1+s^2)^{\frac{3}{2}}}.
$$
Canceling $s$ and passing to  polar coordinates,
$$
d\gamma_{d-2}=
\frac{d}{d-2}\int\limits_{B^{d-2}_2}|p|^2dp=\frac{d}{d-2}\int\limits_{S^{d-3}}d\omega\int\limits_0^1r^{2+d-3}dr=
\frac{\omega(S^{d-3})}{d-2}=\kappa_{d-2},
$$
where $\omega(S^{d-3})$ is the surface area of $S^{d-3}$, we have (\ref{i1}).

Now we prove the converse statement. Fix any $j=2,\dots,d-1$. We rewrite the first equality in (\ref{ij}) as 
$$
\frac{I_j(s)}{ \gamma_{d-2}\, \sqrt{1+s^2}}=\int\limits_{-x(s)}^{y(s)}(f(t)^2-L(s,t)^2)^{\frac{d}{2}}dt
$$
and differentiate both sides with respect to $s$.
Using (\ref{cmass}) and (\ref{i1}), we see that
\begin{equation}\label{rabotav}
\Big(\frac{I_j(s)}{ \sqrt{1+s^2}}\Big)'=\frac{I_j'(s)
	(1+s^2)-sI_j(s)}{ (1+s^2)^{\frac{3}{2}}}=-\frac{const\,\,s}{  (1+s^2)^{\frac{3}{2}}},
\end{equation}
where the second equality above is obtained follows. Using (\ref{defxy}) 
we  differentiate the first equality in (\ref{ij}) to obtain
$$
I_j'(s)(1+s^2)=\gamma_{d-2}\,s\sqrt{1+s^2}\int\limits_{-x(s)}^{y(s)}(f(t)^2-L(s,t)^2)^{\frac{d}{2}}dt-
$$
$$
d\gamma_{d-2}(1+s^2)^{\frac{3}{2}}
\int\limits_{-x(s)}^{y(s)}(f(t)^2-L(s,t)^2)^{\frac{d-2}{2}}(st+h(s))(t+h'(s))dt.
$$
Adding and subtracting $sh'(s)$ in the second parentheses under the second integral  and using (\ref{cmass}), the fact that
 $d\gamma_{d-2}=\kappa_{d-2}$ and the second equality in (\ref{i1}),
we have $$
I_j'(s)(1+s^2)-sI_j(s)=sI_j(s)-sI_1-sI_j(s)=-sI_1=-const\,s.
$$
  This gives the second equality in  (\ref{rabotav}), i.e., 
$$
I_j'(s)-\frac{s}{1+s^2}I_j(s)+const\frac{s}{1+s^2}=0.
$$
Solving  this  linear ODE with an integrating factor $\frac{1}{\sqrt{1+s^2}}$, we have
$$
I_j(s)=\sqrt{1+s^2}\Big(\frac{const}{\sqrt{1+s^2}} +c_1\Big)=const+c_1\,\sqrt{1+s^2}
$$
with some constant $c_1$. Since $I_j$ is bounded on $[s_o,\infty)$, $c_1=0$, and we obtain the converse part of the lemma.
\ep

\medskip

Let 
\begin{equation}\label{ball11}
f_o(t)=\sqrt{1-t^2},\qquad L_o(s,t)=st,\qquad x_o(s)=y_o(s)=\frac{1}{\sqrt{1+s^2}},
\end{equation}
 where $f_o$  describes  the boundary of the unit Euclidean ball, $L_o$ corresponds to the linear subspace  passing through the origin with $h\equiv 0$, 
and  $x_o$, $y_o$ are the first coordinates of the points of intersection of $\pm f$ and $L_o$.
   Our goal is to prove the following proposition.

\bprop\label{hob1}
Let $n=\frac{d}{2}$. A body 
$K_f$ floats in equilibrium in every orientation at the level $\frac{\textnormal{vol}_d(K)}{2}$, provided  for all $s>0$,
\begin{equation}\label{ijn}
\int\limits_{-x(s)}^{y(s)}(f(t)^2-L(s,t)^2)^{n}dt=\int\limits_{-x_o(s)}^{y_o(s)}(f_o(t)^2-L_o(s,t)^2)^{n}dt=\frac{const}{\sqrt{1+s^2}},
\end{equation}
\begin{equation}\label{cmassn}
\int\limits_{-x(s)}^{y(s)}(f(t)^2-L(s,t)^2)^{n-1}\frac{\partial L(s,t)}{\partial s}\,dt=0.
\end{equation}
\eprop
We remark that \eqref{ijn}
and \eqref{cmassn} 
are similar to equations (4) and (5) from \cite{NRZ}.

\bp
Observe that $H(L_0)$ divides $K_f$ into two parts of equal volume. Also, (\ref{cmassn})   is  the same as (\ref{cmass}) of Lemma \ref{FL}. Thus, by
Lemma \ref{charp2} and Lemma  \ref{FL}
the characteristic points of  the family of   hyperplanes $\{H(L_s)$, $s\in [0,\infty)\}$, are exactly the centers of mass 
of the sections $K\cap H(L_s)$. Hence, 
we can apply 
the converse part of Theorem \ref{olovzhal} to conclude  that they 
are the cutting hyperplanes at the level $\frac{\textrm{vol}_d(K)}{2}$. 

On the other hand, observing that  conditions 
 (\ref{ijn}), (\ref{cmassn})  are the same as (\ref{ij}) and (\ref{cmass}), by Lemma \ref{cl1} condition (\ref{i1})
 also holds. 
 Therefore,  by 
Lemma \ref{work1h}, the cutting sections   have equal  moments of inertia  for all 
$(d-2)$-dimensional planes  passing through the centers of mass of these sections.
 By  Remark \ref{rv},  all conditions of
the converse part of Theorem \ref{Fedja1}
are satisfied, and the proposition  follows.
\ep

 In order to construct a counterexample, 
we will choose  the {\it perturbation function} $h$ with the properties described at the beginning of this section. 
The convex body corresponding to any such function will  automatically  be asymmetric since not all its  sections dividing the volume in half will pass through a single point.

\section{The case of even $d\ge 4$}\label{Even}

Note that in this case $n =\frac{d}{2}\in{\mathbb N}$. 	Our argument is very similar to the one  in Section 3 of \cite{NRZ}. Our body  $K_f$ will be  a {\it local} perturbation of 
the Euclidean ball, i.e., 
 the resulting function $f(t)$ will be equal to $\sqrt{1-t^2}$ everywhere on $[-1,1]$ except  $[-\frac{1}{\sqrt{1+(1-2\tau)^2}},-\frac{1}{\sqrt{1+(1-\tau)^2}}]\cup  [\frac{1}{\sqrt{1+(1-\tau)^2}},\frac{1}{\sqrt{1+(1-2\tau)^2}}]  $ for some small $\tau>0$.


Equations (\ref{defxy}) 
show that to define $f$, it is enough to define two decreasing functions $x(s)$, $y(s)$ on $[0, +\infty)$. 
 Our  functions $x(s)$ and $y(s)$ will coincide with $x_o$ and $y_o$  for all $s\notin[1-2\tau,1-\tau]$, where $x_o$, $y_o$ are defined by (\ref{ball11}). Since the curvature of the semicircle is strictly positive, the resulting function $f$ will be strictly concave if $x$ and $y$ are close to $x_o$ and $y_o$ in $C^2$.
	
	We will make our construction in several steps. First, we {\it define} $x=x_o$, $y=y_o$ on $[1,\infty)$. Second, we will express    equations (\ref{ijn}), (\ref{cmassn})  purely in terms of $x$ and $y$ (see  (\ref{hah1}) and (\ref{hah2}) below). Then we
	will use these new equations
	to extend the functions $x$ and $y$ to $[1-3\tau,1]$. We will be able to do it if $\tau$ and $h$ are  sufficiently small. Moreover, the extensions will coincide with $x_o$ and $y_o$ on $[1-\tau,1]$ and will be close to $x_o$ and $y_o$ up to  two derivatives on $[1-3\tau,1-\tau]$. Then, we will show that  our extensions automatically coincide with $x_o$ and $y_o$ on $[1-3\tau,1-2\tau]$ as well. This will allow us to put $x=x_o$, $y=y_o$ on the remaining interval $[0, 1-3\tau]$ and get a nice smooth function. Finally, we will show that equations (\ref{ijn}), (\ref{cmassn}) will be satisfied up to $s=0$, thus finishing the proof.
	
	\bigskip
	
	{\bf Step 1}.
	We put   $x=x_o$, $y=y_o$ on $[1,\infty)$.
	
	\medskip

	{\bf Step 2}. To construct $x$, $y$ on $[1-3\tau,1]$,   we will make  some technical preparations.
	First,  we will  differentiate equations (\ref{ijn}), (\ref{cmassn})  a few times
	to obtain a   system of {\it four} integral  equations with  {\it four} unknown functions $x$, $y$, $x'$, $y'$.
	Next, we will apply Lemma 8 and Remark 2 from \cite[pp. 63-66]{NRZ} to show
	that there exists a  solution $x$,  $y$, $x'$, $y'$ of the constructed system of integral equations on $[1-3\tau,1]$,  which coincides with $x_o$, $y_o$, $\frac{dx_o}{ds}$, $\frac{dy_o}{ds}$ on $[1-\tau,1]$.
	Finally, we will prove that the $x$ and $y$ components of that solution  give a solution of (\ref{ijn}), (\ref{cmassn})  with $f$ defined by (\ref{defxy}).

	Differentiating  equation  (\ref{ijn}) $n+1$ times and   equation   (\ref{cmassn}) $n$ times with respect to $s$ and using (\ref{defxy}),
	we obtain
	$$
	(-2)^nn!\Big[
	\Big(\Big(L\frac{\partial L}{\partial s}\Big)\Big|_{(s,-x(s))}\Big)^n\frac{dx}{ds}(s)\,+\Big(\Big(L\frac{\partial L}{\partial s}\Big)\Big|_{(s,y(s))}\Big)^n\frac{dy}{ds}(s)\Big]\,+
	$$
	\begin{equation}\label{feq}
	\int\limits_{-x(s)}^{y(s)}\Big(\frac{\partial}{\partial s}\Big)^{n+1}\Big((f(t)^2-L(s, t)^2)^n\Big)dt\,=\Big(\frac{d}{ds}\Big)^{n+1}\Big(\frac{const}{\sqrt{1+s^2}}\Big),
	\end{equation}
	and
	$$
	(-2)^{n-1}(n-1)!\Big[\Big(\Big(L\frac{\partial L}{\partial s}\Big)^{n-1}\frac{\partial L}{\partial s}\Big)\Big|_{(s,-x(s))}\frac{dx}{ds}(s)\,+\qquad\qquad\qquad
	$$
	$$
\qquad\qquad\qquad\qquad\qquad\qquad\qquad\qquad	\Big(\Big(L\frac{\partial L}{\partial s}\Big)^{n-1}\frac{\partial L}{\partial s}\Big)\Big|_{(s,y(s))}\frac{dy}{ds}(s)\Big]\,+
	$$
	\begin{equation}\label{seq}
	\int\limits_{-x(s)}^{y(s)}\Big(\frac{\partial}{\partial s}\Big)^{n}\Big((f(t)^2-L(s, t)^2)^{n-1}\frac{\partial L}{\partial s}(s,t)\Big)dt=0.
	\end{equation}
	
	When $s\le 1$, the integral term $I$  in  (\ref{feq}) can be split as
	$$
	I=\int\limits_{-x(s)}^{y(s)}\Big(\frac{\partial}{\partial s}\Big)^{n+1}\Big((f(t)^2-L(s, t)^2)^n\Big)dt=
	$$
	$$
	\Big(\int\limits_{-x(s)}^{-x_o(1)}\,+\,\int\limits_{y_o(1)}^{y(s)}\Big)\,\Big(\frac{\partial}{\partial s}\Big)^{n+1}\Big((f(t)^2-L(s, t)^2)^n              \Big)dt\,+\,\Xi_1(s),
	$$
	where
	$$
	\Xi_1(s)=\int\limits^{y_o(1)}_{-x_o(1)}\Big(\frac{\partial}{\partial s}\Big)^{n+1}\Big((f_o(t)^2-L(s, t)^2)^n\Big)dt.
	$$
	Making the change of variables $t=-x(\sigma)$ in the integral $\int_{-x(s)}^{-x_o(1)}$ and $t=y(\sigma)$  in the integral $\int_{y_o(1)}^{y(s)}$ and using (\ref{defxy}), we obtain
	$$
	I=-\int\limits_s^{1}\Big(\frac{\partial}{\partial s}\Big)^{n+1}\Big(L(\sigma,-x(\sigma))^2-L(s,-x(\sigma))^2\Big)^n\,\frac{dx}{ds}(\sigma)d\sigma\,-
	$$
	$$
	\int\limits_s^{1}\Big(\frac{\partial}{\partial s}\Big)^{n+1}\Big(L(\sigma,y(\sigma))^2-L(s,y(\sigma))^2\Big)^n\,\frac{dy}{ds}(\sigma)d\sigma\,+\,\Xi_1(s).
	$$
	Similarly,  we have
	$$
	\int\limits_{-x(s)}^{y(s)}\Big(\frac{\partial}{\partial s}\Big)^{n}\Big((f(t)^2-L(s, t)^2)^{n-1}\frac{\partial L}{\partial s}(s,t)\Big)dt=
	$$
	$$
	-\int\limits_s^{1}\Big(\frac{\partial}{\partial s}\Big)^{n}\Big(\Big(L(\sigma,-x(\sigma))^2-L(s,-x(\sigma))^2\Big)^{n-1}\frac{\partial L}{\partial s}(s,-x(\sigma))\Big)\,\frac{dx}{ds}(\sigma)d\sigma\,-
	$$
	$$
	\int\limits_s^{1}\Big(\frac{\partial}{\partial s}\Big)^{n}\Big(\Big(L(\sigma,y(\sigma))^2-L(s,y(\sigma))^2\Big)^{n-1}\frac{\partial L}{\partial s}(s,y(\sigma))\Big)\,\frac{dy}{ds}(\sigma)d\sigma\,+\,\Xi_2(s),
	$$
	where
	$$
	\Xi_2(s)=\int\limits^{y_o(1)}_{-x_o(1)}\Big( \frac{\partial}{\partial s}\Big)^{n}\Big((f_o(t)^2-L(s, t)^2)^{n-1}\frac{\partial L}{\partial s}(s,t)                 \Big)dt.
	$$
	To reduce the resulting system of integro-differential equations to a pure system of integral equations we add two independent unknown functions $x'$, $y'$ and
	two new relations:
	$$
	  x(s)=-\int\limits_s^{1}x'(\sigma)d\sigma +x_o(1), \qquad y(s)=-\int\limits_s^{1}y'(\sigma)d\sigma +y_o(1).
	$$
	We rewrite   our  equations (\ref{feq}), (\ref{seq}) as follows:
	\begin{equation}\label{hah1}
	(-2)^nn!\Big[\Big(\Big(L\frac{\partial L}{\partial s}\Big)\Big|_{(s,-x(s))}\Big)^nx'(s)\,+
	\Big(\Big(L\frac{\partial L}{\partial s}\Big)\Big|_{(s,y(s))}\Big)^ny'(s)\Big]\,-
	\end{equation}
	$$
	\int\limits_s^{1}\Big(\frac{\partial}{\partial s}\Big)^{n+1}\Big(L(\sigma,-x(\sigma))^2-L(s,-x(\sigma))^2\Big)^n\,x'(\sigma)d\sigma\,
	-
	$$
	$$
	\int\limits_s^{1}\Big(\frac{\partial}{\partial s}\Big)^{n+1}\Big(L(\sigma,y(\sigma))^2-L(s,y(\sigma))^2\Big)^n\,y'(\sigma)d\sigma\,+\,\Xi_1(s)=
	$$
	$$
	\Big(\frac{d}{ds}\Big)^{n+1}\Big(\frac{const}{\sqrt{1+s^2}}\Big),
	$$
	and
	\begin{equation}\label{hah2}
	(-2)^{n-1}(n-1)!\Big[\Big(\Big(L\frac{\partial L}{\partial s}\Big)^{n-1}\frac{\partial L}{\partial s}\Big)\Big|_{(s,-x(s))}x'(s)\,+\qquad\qquad\qquad
		\end{equation}
	$$
\qquad\qquad\qquad\qquad\qquad\qquad\qquad\qquad	\Big(\Big(L\frac{\partial L}{\partial s}\Big)^{n-1}\frac{\partial L}{\partial s}\Big)\Big|_{(s,y(s))}y'(s)\Big]
	\,-
$$
	$$
	\int\limits_s^{1}\Big(\frac{\partial}{\partial s}\Big)^{n}\Big(\Big(L(\sigma,-x(\sigma))^2-L(s,-x(\sigma))^2\Big)^{n-1}\frac{\partial L}{\partial s}(s,-x(\sigma))\Big)\,x'(\sigma)d\sigma\,-
	$$
	$$
	\int\limits_s^{1}\Big(\frac{\partial }{\partial s}\Big)^{n}\Big(\Big(L(\sigma,y(\sigma))^2-L(s,y(\sigma))^2\Big)^{n-1}\frac{\partial L}{\partial s}(s,y(\sigma))\Big)\,y'(\sigma)d\sigma+	\Xi_2(s)=0.
	$$
	
	Now we  rewrite our system  in the  form
	\begin{equation}\label{hah41}
	\mathbf{G}(s,Z(s))=\int\limits_s^{1}\mathbf{\Theta}(s,\sigma, Z(\sigma))d\sigma+\mathbf{\Xi}(s).
	\end{equation}
	Here
	$$
	Z=\left(
	\begin{array}{cccc}
	x\\
	y\\
	x'\\
	y'
	\end{array}
	\right),
	$$
	$$
	\mathbf{G}(s,Z)=\left(
	\begin{array}{cccc}
	x\\
	y\\[7pt]
	(-2)^nn!\Big[\Big(L\frac{\partial L}{\partial s}\Big|_{(s,-x)}\Big)^nx'\,+
	\Big(L\frac{\partial L}{\partial s}\Big|_{(s,y)}\Big)^ny'\Big]\\[15pt]
	(-2)^{n-1}(n-1)!\Big[\Big(\Big(L\frac{\partial L}{\partial s}\Big)^{n-1}\frac{\partial L}{\partial s}\Big)\Big|_{(s,-x)}x'\,+
	\Big(\Big(L\frac{\partial L}{\partial s}\Big)^{n-1}\frac{\partial L}{\partial s}\Big)\Big|_{(s,y)}y'\Big]\\
	\end{array}
	\right),
	$$
	$$
	\mathbf{\Theta}(s,\sigma, Z)=
	-\left(
	\begin{array}{cccc}
	x'\\
	y'\\
	\Theta_{1}\\
	\Theta_{2}
	\end{array}
	\right),
	$$
	where
	$$
	\Theta_{1}=-\Big(\frac{\partial}{\partial s}\Big)^{n+1}\Big(L(\sigma,-x)^2-L(s,-x)^2\Big)^n\,x'\,
	-
	\Big(\frac{\partial}{\partial s}\Big)^{n+1}\Big(L(\sigma,y)^2-L(s,y)^2\Big)^n\,y'\,,
	$$
	$$
	\Theta_{2}=
	-\Big(\frac{\partial}{\partial s}\Big)^{n}\Big(\Big(L(\sigma,-x)^2-L(s,-x)^2\Big)^{n-1}\frac{\partial L}{\partial s}(s,-x)\Big)\,x'\,-\qquad\qquad
	$$
	$$
	\qquad\qquad\qquad\qquad\qquad\qquad\qquad\Big(\frac{\partial}{\partial s}\Big)^{n}\Big(\Big(L(\sigma,y)^2-L(s,y)^2\Big)^{n-1}\frac{\partial L}{\partial s}(s,y)\Big)\,y',
	$$
	and
	$$
	\mathbf{\Xi}(s)=\left(
	\begin{array}{cccc}
	x_o(1)\\
	y_o(1)\\
	-\Xi_1(s)+\Big(\frac{d}{ds}\Big)^{n+1}\Big(\frac{const}{\sqrt{1+s^2}}\Big)\\
	-\Xi_2(s)
	\end{array}
	\right).
	$$
	Note that ${\mathbf G}$, ${\mathbf \Theta}$, ${\mathbf \Xi}$ are well defined and infinitely smooth for all $s,\sigma \in(0,1]$ and $Z\in \R^4$.
	Observe also that
	$$
	D_Z{\mathbf G}\Big|_{(s,Z)}=
	\left(
	\begin{array}{cc}
	\mathbf{I}&0\\
	\mathbf{*}&\mathbf{A}\\
	\end{array}
	\right),
	$$
	where
	$$
	{\mathbf I}=\left(
	\begin{array}{ccc}
	1 & 0\\
	0&
	1
	\end{array}
	\right),\qquad{\mathbf A}={\mathbf A(s,x,y)}=
	$$
	$$
	\left(
	\begin{array}{ccc}
	(-2)^nn!\Big(\Big(L\frac{\partial L}{\partial s}\Big)\Big|_{(s,-x)}\Big)^n\qquad& (-2)^nn!\Big(\Big(L\frac{\partial L}{\partial s}\Big)\Big|_{(s,y)}\Big)^n\\
	(-2)^{n-1}(n-1)!\Big(\Big(L\frac{\partial L}{\partial s}\Big)^{n-1}\frac{\partial L}{\partial s} \Big)\Big|_{(s,-x)}\qquad&
	(-2)^{n-1}(n-1)!\Big(\Big(L\frac{\partial L}{\partial s}\Big)^{n-1}\frac{\partial L}{\partial s}\Big)\Big|_{(s,y)}
	\end{array}
	\right).
	$$
	
	\smallskip

	The function
	$$
	Z_o(s)=\left(
	\begin{array}{cccc}
	x_o(s)\\
	y_o(s)\\[3pt]
	\frac{dx_o}{ds}(s)\\[3pt]
	\frac{dy_o}{ds}(s)
	\end{array}
	\right)
	$$
	solves the system (\ref{hah41}) with ${\mathbf G}$, ${\mathbf \Theta}$, ${\mathbf \Xi}$ corresponding to $h\equiv 0$ (we will denote them by
	${\mathbf G_o}$, ${\mathbf \Theta_o}$, ${\mathbf \Xi_o}$) on $[\frac{1}{2},1]$.

	We claim that
	\begin{equation}\label{detH}
	\det\,\Big(D_Z{\mathbf G_o}\Big|_{(s,Z_o(s))}  \Big)=\det({\mathbf A_o(s,x_o(s),y_o(s))})\neq 0\qquad \forall s\in(0,1].
	\end{equation}
	Indeed,  since  the  matrix ${\mathbf A_o}(s,x_o(s),y_o(s))$
	is of  the form
	$$
	\left(
	\begin{array}{ccc}
	(-2)^nn!(sx_o(s))^n\qquad& (-2)^nn!(sy_o(s))^n\\
	(-2)^{n-1}(n-1)!(sx_o(s))^{n-1}(-x_o(s))\qquad&
	(-2)^{n-1}(n-1)!(sy_o(s))^{n-1}y_o(s)
	\end{array}
	\right),
	$$
its sign pattern is
	$$
	\left(
	\begin{array}{ccc}
	+ & +\\
	+ &  -
	\end{array}
	\right),
	\qquad \textrm{when}\,\, n\,\,\textrm{ is even},\quad \textrm{and }
	\qquad
	\left(\begin{array}{ccc}
	- & -\\
	-&
	+
	\end{array}
	\right),\quad \textrm{when}\, n\,\,\textrm{ is odd}.
	$$
	Thus, (\ref{detH}) follows.
	In particular,
	$$
	\det\,\Big(D_Z{\mathbf G_o}\Big|_{(1,Z_o(1))}  \Big)\neq 0.
	$$

	Lemma 8 from \cite[p. 63]{NRZ} then  implies that we can choose some small $\tau>0$ and, for any fixed $k\in {\mathbb N}$,  construct a solution $Z(s)$ of (\ref{hah41}) which is 
	$C^k$-close to $Z_o(s)$  on $[1-3\tau,1]$, whenever
	${\mathbf G}$, ${\mathbf \Theta}$, ${\mathbf \Xi}$ are sufficiently close to ${\mathbf G_o}$, ${\mathbf \Theta_o}$, ${\mathbf \Xi_o}$ in $C^k$ on  certain compact sets. Since
	$\mathbf{G}$, $\mathbf{\Theta}$, $\mathbf{\Xi}$ and their derivatives are some explicit (integrals of) polynomials in $Z$, $s$, $\sigma$, $h(s)$, and the derivatives of $h(s)$,
	this closeness condition will hold if $h$ and sufficiently many of its derivatives are close enough to zero.
	Moreover, since $h$ vanishes on $[1-\tau,1]$, the assumptions of Remark 2 from \cite[p. 66]{NRZ} are satisfied and we have $Z(s)=Z_o(s)$ on $[1-\tau,1]$.

	To prove that the $x$ and $y$ components of the solution we found  give a solution of (\ref{ijn}), (\ref{cmassn})  with $f$ defined by (\ref{defxy}), we
	consider the functions
	$$
	F(s):=
	\int\limits_{-x(s)}^{y(s)}\Big(f(s,t)^2-L(s,t)^2\Big)^ndt-\frac{const}{\sqrt{1+s^2}},
	$$
	$$
	H(s):=  \int\limits_{-x(s)}^{y(s)}\Big(f(s,t)^2-L(s,t)^2\Big)^{n-1}\frac{\partial L}{\partial s}(s,t)dt.
	$$
	Since equations (\ref{hah1}) and  (\ref{hah2}) of our system (\ref{hah41}) were obtained by the  differentiation of equations (\ref{ijn}),  (\ref{cmassn}),  we have
	$$
	\Big(\frac{d}{ds}\Big)^{n+1}F(s)=0,\qquad \Big(\frac{d}{ds}\Big)^nH(s)=0
	$$
	on $[1-3\tau, 1]$. Hence, $F$ and $H$ are polynomials
	on $[1-3\tau, 1]$.
	Since $h(s)=0$, $x(s)=x_o(s)$, $y(s)=y_o(s)$  on  $[1-\tau, 1]$,  $F$ and $H$  vanish on $[1-\tau, 1]$ and, therefore, identically.
	Thus, we conclude
	that  the $x$ and $y$ components of    the solutions
	of (\ref{hah1}), (\ref{hah2}) solve  (\ref{ijn}),  (\ref{cmassn})  on $(1-3\tau,1]$.
	Step 2 is completed.

	\medskip

	{\bf Step 3}.
	We claim that $x=x_o$, $y=y_o$ on $[1-3\tau, 1-2\tau]$, i.e., the {\it perturbed solution} returns to the semicircle.
	Since  $h$ is supported on $[1-2\tau,1-\tau]$,  we have $L=L_o=st$ and
	$\frac{\partial}{\partial s}L(s,t)=t$ for $s\in [1-3\tau,1-2\tau]$. It follows that every time we differentiate equation   (\ref{ijn})  (with respect to $s$)
	we can divide the result by $s$ to obtain
	\begin{equation}\label{evenmom}
	\int\limits_{-x(s)}^{y(s)}(f(t)^2-L_o(s,t)^2)^{n-k}t^{2k}dt=
	\int\limits_{-x_o(s)}^{y_o(s)}(f_o(t)^2-L_o(s,t)^2)^{n-k}t^{2k}dt,
	\end{equation}
for $k\le n$.	If we take $k=n$ in (\ref{evenmom}), we get
	\begin{equation}\label{evenk1}
	\int\limits_{-x(s)}^{y(s)}t^{2n}dt=
	\int\limits_{-x_o(s)}^{y_o(s)}t^{2n}dt.
	\end{equation}
	Similarly, for $k\le n-1$,   equation (\ref{cmassn}) implies that
	\begin{equation}\label{oddmom}
	\int\limits_{-x(s)}^{y(s)}(f(t)^2-L_o(s,t)^2)^{n-1-k}t^{2k+1}dt=\quad\qquad\qquad\qquad\qquad\qquad\qquad
		\end{equation}
		$$
	\qquad\qquad\qquad\qquad\qquad\qquad\int\limits_{-x(s)}^{y(s)}(f(t)^2-L_o(s,t)^2)^{n-1-k}t^{2k+1}dt=0.
$$
	Putting $k=n-1$ in (\ref{oddmom}), we get
	\begin{equation}\label{oddk}
	\int\limits_{-x(s)}^{y(s)}t^{2n-1}dt=0=
	\int\limits_{-x_o(s)}^{y_o(s)}t^{2n-1}dt.
	\end{equation}
	Equation (\ref{oddk}) yields $x(s)=y(s)$, and the symmetry (with respect to $0$) of the intervals $(-x_o(s),y_o(s))$, $(-x(s),y(s))$,
	together with (\ref{evenk1}), yield $(-x_o(s),y_o(s))$ $=(-x(s),y(s))$ for all $s\in [1-3\tau, 1-2\tau]$.
	Step 3 is completed.
	
	\medskip

	{\bf Step 4}.  We put $x=x_o$, $y=y_o$ on $ [0,1-3\tau]$, which will result in a function $f$ defined  on $[-1,1]$ and coinciding with $f_o(t)=\sqrt{1-t^2}$ outside
	small intervals around $\pm\frac{1}{\sqrt{2}}$.
	It remains to check that (\ref{ijn}),  (\ref{cmassn}) are valid for $s\in[0,1-3\tau]$.  We  will prove the validity of  (\ref{cmassn}). The proof
	for  equation (\ref{ijn})  is similar and can be found in \cite[p. 53]{NRZ}.

Since $h\equiv 0$ away from $(1-2\tau,1-\tau)$, we have $L(s,t)=st$ for $s\in[0,1-3\tau]$,  so we need to check that
	$$
	\int\limits_{-x(s)}^{y(s)}(f(t)^2-(st)^2)^{n-1}\,tdt=\int\limits_{-x(s)}^{y(s)}(f_o(t)^2-(st)^2)^{n-1}\,tdt,\qquad \forall s\in[0,1-3\tau].
	$$
	Recall that $x=x_o$ and $y=y_o$ everywhere on this interval, so we can write $x$ and $y$ instead of $x_o$ and $y_o$ on the right-hand side.

	Using the binomial formula, we see that it suffices to check that
	\begin{equation}\label{devenm}
	\int\limits_{-x(s)}^{y(s)}f(t)^{2j}\,t^{2(n-1-j)+1}dt=\int\limits_{-x(s)}^{y(s)}f_o(t)^{2j}\,t^{2(n-1-j)+1}dt, 
	\end{equation}
$\forall j=1,\dots,n-1$  and	$s\in[0,1-3\tau]$.
	Since $f\equiv f_o$  outside $[-x(1-3\tau),y(1-3\tau)]$,
	splitting the integrals in (\ref{devenm}) into three parts
 with ranges $[-x(s),-x(1-3\tau)]$, $[-x(1-3\tau),y(1-3\tau)]$, $[y(1-3\tau),y(s)]$,
	 it is enough to check (\ref{devenm}) on the middle interval $[-x(1-3\tau),y(1-3\tau)]$.
	
	To this end, we first  take $s=1-3\tau$, $k=n-2$ in (\ref{oddmom}) and conclude that
	\begin{equation}\label{evencase2}
	\int\limits_{-x(1-3\tau)}^{y(1-3\tau)}f(t)^{2}\,t^{2n-3}dt=\int\limits_{-x(1-3\tau)}^{y(1-3\tau)}f_o(t)^{2}\,t^{2n-3}dt,
	\end{equation}
	which is (\ref{devenm}) for $j=1$ and $s=1-3\tau$.
	Now we go ``one step up",  by taking $s=1-3\tau$, $k=n-3$ in (\ref{oddmom}), to get
	$$
	\int\limits_{-x(1-3\tau)}^{y(1-3\tau)}(f(t)^{2}-(st)^2)^2t^{2n-5}dt=\int\limits_{-x(1-3\tau)}^{y(1-3\tau)}(f_o(t)^{2}-(st)^2)^2t^{2n-5}dt.
	$$
	The last equality  together with (\ref{evencase2}) yield
	$$
	\int\limits_{-x(1-3\tau)}^{y(1-3\tau)}f(t)^{4}\,t^{2n-5}dt=\int\limits_{-x(1-3\tau)}^{y(1-3\tau)}f_o(t)^{4}t^{2n-5}dt,
	$$
	which is (\ref{devenm}) for $j=2$ and $s=1-3\tau$. Proceeding in a similar way we get (\ref{devenm}) for $j=1,\dots,n-1$ and $s=1-3\tau$.
	This finishes the proof of Theorem \ref{bitsya} in even dimensions. $\qquad\,\,\,\qquad\qquad\qquad\qquad\qquad\qquad\qquad\qquad\qquad\qquad\square$

\section{The case of odd $d\ge 3$}\label{Odd}


Note that  $n=q+\frac{1}{2}$, $q\in{\mathbb N}$.
Then (\ref{ijn}) and (\ref{cmassn}) take the form
\begin{equation}\label{dodd+}
\int\limits_{-x(s)}^{y(s)}(f(t)^2-L(s,t)^2)^{q+\frac{1}{2}}dt=\int\limits_{-x_o(s)}^{y_o(s)}(f_o(t)^2-L_o(s,t)^2)^{q+\frac{1}{2}}dt=\frac{const}{\sqrt{1+s^2}},
\end{equation}
\begin{equation}\label{dodd-}
\int\limits_{-x(s)}^{y(s)}(f(t)^2-L(s,t)^2)^{q-\frac{1}{2}}\frac{\partial L}{\partial s}(s,t)dt =0,
\end{equation}
where $f_o$, $L_o$,   $y_o$,  and $x_o$ are defined by (\ref{ball11}).

Our argument is  similar to the one  in \cite[Section 4]{NRZ}. 
Our body of revolution $K_f$ will be constructed as a  perturbation of the Euclidean ball. We remark  that in the case of  odd dimensions,  the perturbation will not  be local,  
meaning that  the resulting  function $f(t)$ will be equal to $\sqrt{1-t^2}$   on $\left[-\frac{1}{\sqrt{1+(1-\tau)^2}},\frac{1}{\sqrt{1+(1-\tau)^2   }}\right]$ for some small $\tau>0$.

We will make our construction in several steps  corresponding to the slope ranges $s\in [1,\infty)$, $s\in[1-3\tau,1]$, and $s\in (0,1-3\tau]$.
We will use different   ways to describe the boundary of $K_f$ within those ranges. We will {\it define}
$f(t)=f_o(t)$ for  $t\in\left[-\frac{1}{\sqrt{2}},\frac{1}{\sqrt{2}}\right]$.
We will differentiate    (\ref{dodd+}), (\ref{dodd-}) and rewrite  the resulting equations in terms of $x$ and $y$, to extend
$x$ and $y$ to $[1-3\tau, 1]$ like we did in the even case. As before, $f$ is related to $x$ and $y$ by (\ref{defxy}).
Finally, we will change the point of view and define the remaining part of $f$ in terms of the  functions
$R(\alpha)$ and $r(\alpha)$, related to $f$  by
\begin{equation}\label{f}
f(R(\alpha)\cos\alpha)=R(\alpha)\sin\alpha,\,\, f(-r(\alpha)\cos\alpha)=r(\alpha)\sin\alpha,\,\, \alpha\in[0,\tfrac{\pi}{2}].
\end{equation}
Note that the {\it radial function }
$\rho_K(w)=\sup\{t>0:\,\,tw\in K\}$
of the resulting body $K$ satisfies
\begin{equation}\label{tryuk}
\rho_K(w)= \begin{cases}
R(\alpha)\qquad &\text{if \,}w_1>0,\\
r(\alpha) & \text{if \,}w_1<0,
\end{cases}
\end{equation}
where $w=(w_1,\dots, w_d)\in S^{d-1}$ and $\alpha\in[0,\frac{\pi}{2}]$, $\cos\alpha=|w_1|$.

\bigskip

{\bf Step 1}.
We put $x=x_o$, $y=y_o$ on $[1,\infty)$, which is equivalent to putting $f(t)=\sqrt{1-t^2}$ for $t\in[-\frac{1}{\sqrt{2}},\frac{1}{\sqrt{2}}]$.

\medskip

{\bf Step 2}.
Differentiating  equation  (\ref{dodd+}) $q+1$ times, we obtain
\begin{equation}\label{ass+}
\Big(\frac{\partial}{\partial s}  \Big)^{q+1}\int\limits_{-x(s)}^{y(s)}(f(t)^2-L(s,t)^2)^{q+\frac{1}{2}}dt=
\end{equation}
$$
\Big(\int\limits_{-x(s)}^{-x_o(1)}\,+\,\int\limits_{y_o(1)}^{y(s)}\Big)\Big(\frac{\partial}{\partial s}  \Big)^{q+1}\Big((f(t)^2-L(s,t)^2)^{q+\frac{1}{2}}\Big)dt\,
+E_1(s)=
$$
$$
\Big(\frac{d}{ds}  \Big)^{q+1}\frac{const}{\sqrt{1+s^2}},
$$
where
$$
E_1(s)=\int\limits_{-x_o(1)}^{y_o(1)}\Big(\frac{\partial}{\partial s}  \Big)^{q+1}\Big((f_o(t)^2-L(s,t)^2)^{q+\frac{1}{2}}\Big)dt.
$$
Note that, unlike  the function $\Xi_1$ in the even-dimensional case,  $E_1$ is well defined only for $s\le 1$ and only if $ \|h\|_{C^1}$  is much smaller than $1$. Also,  even with these assumptions, $E_1(s)$ is $C^{\infty}$ on $[0, 1)$ but not at $1$, where it is merely continuous.

Observe  that
$$
\Big(\frac{\partial}{\partial s}  \Big)^{q+1}\Big((f(t)^2-L(s,t)^2)^{q+\frac{1}{2}}\Big)=\frac{J_1(s,t, f(t))}{\sqrt{f^2(t)-L^2(t)}},
$$
where $J_1(s,t,f)$ is some polynomial expression in $s$, $t$, $f$, $h(s)$, and the derivatives of $h$ at $s$.

Making the change of variables $t=-x(\sigma)$ in the  integral $\int_{-x(s)}^{-x_o(1)}$, and $t=y(\sigma)$  in the integral $\int_{y_o(1)}^{y(s)}$ and using (\ref{defxy}), we can rewrite the sum of integrals on the left-hand side of (\ref{ass+}) as
$$
-\,\int\limits_s^{1}\!\Big[\frac{J_1(s,-x(\sigma), L(\sigma,-x(\sigma)))}{\sqrt{L(\sigma, -x(\sigma))^2-L(s, -x(\sigma))^2}}\frac{dx}{ds}(\sigma)+\frac{J_1(s,y(\sigma), L(\sigma, y(\sigma)))}{\sqrt{L(\sigma, y(\sigma))^2-L(s, y(\sigma))^2}}\frac{dy}{ds}(\sigma)\!\Big]\!
d\sigma.
$$
Now write
$$
L(\sigma,t)^2-L(s, t)^2=(L(\sigma,t)-L(s, t))(L(\sigma,t)+L(s, t)),
$$
and
$$
L(\sigma,t)-L(s, t)=\sigma t+h(\sigma)-st-h(s)=
(\sigma-s)(t+H(s,\sigma)),
$$
where
$$
H(s,\sigma)=\frac{h(\sigma)-h(s)}{\sigma-s}=\int\limits_0^1h'(s+(\sigma-s)\tau)d\tau
$$
is an infinitely smooth function of $s$ and $\sigma$.
Let
$$
K_1(s,\sigma, t)=\frac{J_1(s, t, L(\sigma, t))}{\sqrt{ (t+ H(s,\sigma))(L(\sigma,t)+L(s, t))   }   }.
$$
The function $K_1$ is well defined and infinitely smooth for all $s$, $\sigma$, $t$ satisfying
$ (t+ H(s,\sigma))$$(L(\sigma,t)+L(s, t))>0$. If $\|h\|_{C^1}$ is small enough, this condition is fulfilled whenever  $s$, $\sigma$$\in[\frac{1}{2},1]$ and $|t|>\frac{1}{2}$.

Now we can rewrite equation (\ref{ass+}) in the form
\begin{equation}\label{t1}
-\,\int\limits_s^{1}\Big(K_1(s,\sigma, -x(\sigma))
\,\frac{dx}{ds}(\sigma)\,+\,K_1(s,\sigma, y(\sigma))\,
\frac{dy}{ds}(\sigma)\Big)\,\frac{d\sigma}{\sqrt{\sigma-s}}=
\end{equation}
$$
-E_1(s)+\Big(\frac{d}{ds}  \Big)^{q+1}\frac{const}{\sqrt{1+s^2}}.
$$

Similarly, we can differentiate
(\ref{dodd-}) $q$ times and transform the resulting equation
into
\begin{equation}\label{t2}
-\,\int\limits_s^{1}\Big(K_2(s,\sigma, -x(\sigma))
\,\frac{dx}{ds}(\sigma)\,+\,K_2(s,\sigma, y(\sigma))\,
\frac{dy}{ds}(\sigma)\Big)\,\frac{d\sigma}{\sqrt{\sigma-s}}=
\end{equation}
$$
=-E_2(s),
$$
where $K_2$ is well defined and infinitely smooth in the same range as $K_1$.

The function $E_2$ on the right-hand side  of (\ref{t2})   is given by
$$
E_2(s)=\int\limits_{-x_o(1)}^{y_o(1)}\Big(\frac{\partial}{\partial s}  \Big)^{q}\Big((f_o(t)^2-L(s,t)^2)^{q-\frac{1}{2}}\frac{\partial L}{\partial s}(s,t)\Big)dt,
$$
and everything that we said about $E_1$ applies to $E_2$ as well.

Equations (\ref{t1}) and (\ref{t2}) together can be written in the form
\begin{equation}\label{arik}
\int\limits_s^{1}\frac{K(s,\sigma,z(\sigma),\frac{dz}{ds}(\sigma))}{\sqrt{\sigma-s}}d\sigma=Q(s),
\end{equation}
where, for $z=\left(
\begin{array}{cc}
x\\
y
\end{array}
\right),\,\,z'=\left(
\begin{array}{cc}
x'\,\\
y'\,
\end{array}
\right)\in\R^2$,
$$
K(s,\sigma,z,z')=\,-\,\left(
\begin{array}{cc}
K_1(s,\sigma, -x)
\,x'\,+\,K_1(s,\sigma, y)\,
y'\,\\
K_2(s,\sigma, -x)
\,x'\,+\,K_2(s,\sigma, y)\,
y'\,
\end{array}
\right),
$$
$$
Q(s)=\left(
\begin{array}{cc}
-E_1(s)+\Big(\frac{d}{ds}  \Big)^{q+1}\frac{const}{\sqrt{1+s^2}}\,\\[10pt]
-E_2(s)\,
\end{array}
\right).
$$

By Lemma  8 in \cite[p. 63]{NRZ} with $b=1$,   equation (\ref{arik}) is equivalent to
\begin{equation}\label{dvi}
-G_2(s,s,z,z')+\int\limits_s^{1}\frac{\partial}{\partial s}G_2(s,\sigma, z(\sigma),\frac{dz}{ds}(\sigma))d\sigma=\widetilde{Q}(s),
\end{equation}
where
$$
G_2(s,\sigma,z,z')=\int\limits_0^1\frac{K(s+\tau(\sigma-s),\sigma,z,z')}{\sqrt{\tau(1-\tau)}}d\tau,\qquad   \widetilde{Q}(s)= \frac{d}{ds}\int\limits_s^{1}\frac{Q(s')}{\sqrt{s'-s}}ds'.
$$
Note that
$$
G_2(s,s,z,z')=C\,\cdot\,K(s,s,z,z'),\qquad C=\int\limits_0^1\frac{d\tau}{\sqrt{\tau(1-\tau)}}.
$$
To reduce the resulting system of integro-differential equations to a pure system of integral equations we add two independent unknown functions $x'$, $y'$, let
$z'=\left(
\begin{array}{cc}
x'\\
y'
\end{array}
\right)$, $z_o(s)=\left(
\begin{array}{cc}
x_o(s)\\
y_o(s)
\end{array}
\right)$, and
add
two new relations
$$
z(s)=-\int\limits_s^{1}z'(\sigma)d\sigma +z_o(1).
$$
Together with  (\ref{dvi}), they  lead to the system
\begin{equation}\label{hah44}
\mathbf{G}(s,Z(s))=\int\limits_s^{1}\mathbf{\Theta}(s,\sigma, Z(\sigma))d\sigma+\mathbf{\Xi}(s),\qquad Z=\left(
\begin{array}{cc}
z\\
z'
\end{array}
\right)=
\left(
\begin{array}{cccc}
x\\
y\\
x'\\
y'
\end{array}
\right).
\end{equation}
Here
$$
\mathbf{G}(s,Z)=\left(
\begin{array}{cc}
z\\
-G_2(s,s,z,z')\\
\end{array}
\right),\qquad
\mathbf{\Theta}(s,\sigma, Z)=-\,
\left(
\begin{array}{cc}
z'\\
\frac{\partial}{\partial s}G_2(s,\sigma,z,z')\\
\end{array}
\right),
$$
and
$$
\mathbf{\Xi}(s)=\left(
\begin{array}{cc}
z_o(1)\\
\widetilde{Q}(s)
\end{array}
\right).
$$

In what follows, we will choose
$h$  so that   $\|h\|_{C^1}$ is much smaller than $1$.
In this case,
${\mathbf G}$, ${\mathbf \Theta}$ are well defined and infinitely smooth whenever $s$, $\sigma$ $\in [\frac{1}{2},1]$, $|x|$, $|y|$$>\frac{1}{2}$, $z'\in \R^2$, and $
{\mathbf \Xi}$ is well defined  and infinitely smooth on $[\frac{1}{2},1)$.
Observe also  that
$$
D_Z{\mathbf G}\Big|_{(s,Z(s))}=
\left(
\begin{array}{cc}
\mathbf{I}&0\\
\mathbf{*}&\mathbf{A}\\
\end{array}
\right),
$$
where
$$
{\mathbf I}=\left(
\begin{array}{ccc}
1 & 0\\
0&
1
\end{array}
\right),\qquad
{\mathbf A(s,z)}=C\,\cdot\,{\mathbf E(s,z)},
$$
and
$$
{\mathbf E(s,z)}=
\left(
\begin{array}{cc}
K_1(s,s, -x)\,&\,
K_1(s,s, y)\,
\\
K_2(s,s, -x)
\,&\,\,K_2(s,s, y)\,
\end{array}
\right).
$$

The function
$$
Z_o(s)=\left(
\begin{array}{cc}
z_o(s)\\[3pt]
\frac{dz_o}{ds}(s)
\end{array}
\right)
=
\left(
\begin{array}{cccc}
x_o(s)\\
y_o(s)\\[3pt]
\frac{dx_o}{ds}(s)\\[3pt]
\frac{dy_o}{ds}(s)
\end{array}
\right)
$$
solves the system (\ref{hah44}) with ${\mathbf G}$, ${\mathbf \Theta}$, ${\mathbf \Xi}$ corresponding to $h\equiv 0$ (we will denote them by
${\mathbf G_o}$, ${\mathbf \Theta_o}$, ${\mathbf \Xi_o}$) on $[\frac{1}{2},1]$, say.

We claim that
\begin{equation}\label{detH111}
\det\,\Big(D_Z{\mathbf G_o}\Big|_{(s,Z_o(s))}  \Big)=\det({\mathbf A_o(s,z_o(s))})\neq 0\qquad \textrm{for all}\,\,\,s\in[\tfrac{1}{2},1].
\end{equation}
Indeed,
since $K_{1,2}(s,s,t)$  have the same signs as $J_{1,2}(s,\xi,L(s,t))$ and  since
$$
J_1(s,t, L(s,t))=(2q+1)!!\,\Big(-L(s,t)\,\frac{\partial L}{\partial s}(s,t)\Big)^{q+1},
$$
$$
J_2(s,t, L(s,t))=(2q-1)!!\,\Big(-L(s,t)\frac{\partial L}{\partial s}(s,t)\Big)^{q}\frac{\partial L}{\partial s}(s,t),
$$
we conclude that the matrix ${\mathbf A_o}(s,z_o(s))$ has the same sign pattern as the matrix
$$
\left(
\begin{array}{ccc}
(-1)^{q+1}& (-1)^{q+1}\\
(-1)^{q}(-x_o(s)) &  (-1)^{q}y_o(s)
\end{array}
\right),
$$
i.e.,  the signs in the first row are the same, and the signs in the second one are opposite.

Thus, (\ref{detH111}) follows.
In particular,
$$
\det\,\Big(D_Z{\mathbf G_o}\Big|_{(1,Z_o(1))}  \Big)\neq 0.
$$

Lemma 8 from \cite[p. 63]{NRZ}  then implies  that we can choose some small $\tau>0$ and construct a $C^k$-close to $Z_o(s)$ solution $Z(s)$ of (\ref{hah44}) on $[1-3\tau,1]$ whenever
${\mathbf G}$, ${\mathbf \Theta}$, ${\mathbf \Xi}$ are sufficiently close to ${\mathbf G_o}$, ${\mathbf \Theta_o}$, ${\mathbf \Xi_o}$ in $C^k$ on  certain compact sets. Since
$\mathbf{G}$, $\mathbf{\Theta}$, $\mathbf{\Xi}$ and their derivatives are   (integrals of)  explicit elementary expressions  in $Z$, $s$, $\sigma$, $h(s)$, and the derivatives of $h(s)$,
this closeness condition will hold if $h$ and sufficiently many of its derivatives are close enough to zero.
Moreover, since $h$ vanishes on $[1-\tau,1]$, the assumptions of Remark 2 from \cite[p. 66]{NRZ} are satisfied and we have $Z(s)=Z_o(s)$ on $[1-\tau,1]$.

The $x$ and $y$ components of $Z$ solve the equations obtained by differentiating (\ref{dodd+}) and (\ref{dodd-}). The passage to  (\ref{dodd+}),  (\ref{dodd-}) is now exactly the same as in the even case.

\medskip

{\bf Step 3}.  From now on,  we change the point of view and switch to  the functions $R(\alpha)$ and $r(\alpha)$, $\alpha\in (0,\frac{\pi}{2})$, related
to $f$ by (\ref{f}). The functions $x$ and $y$, which  we have already constructed, implicitly define  $C^{\infty}$-functions $R_h(\alpha)$ and $r_h(\alpha)$ for all $\alpha$ with $\tan\alpha>1-3\tau$.

Instead of parameterizing hyperplanes by the slopes $s$ of the corresponding linear functions, we will
parameterize them by the angles $\beta$ they make with the $x_1$-axis, where $\beta$ is related to $s$ by $\tan\beta=s$.

Our next task will be to derive the equations that will ensure that all {\it central} sections corresponding to angles $\beta$ with $\tan\beta<1-2\tau$ are the cutting sections with  equal  moments with respect to any $(d-2)$-dimensional subspace passing through the origin. We will also ensure that the origin is  the center of mass of these sections.
Note that the sections  are already defined and satisfy these properties when
$\tan \beta\in (1-3\tau,1-2\tau)$.

It will be convenient to  rewrite conditions (\ref{cmass}), (\ref{i1}) and (\ref{ij}) 
in terms of the spherical Radon transform (see \cite[pp. 427-436]{Ga}), defined as
$$
{\mathcal R}f(\xi)=\int\limits_{S^{d-1}\cap\xi^{\perp}}f(w)dw,\qquad f\in C(S^{d-1}),\qquad \xi\in S^{d-1}.
$$

We will use the following  proposition.
\bprop\label{workHARD1}
Let $K$ be a convex body of revolution about the $x_1$-axis containing  the origin  in its interior  and let $\xi=(\pm\sin\alpha,0,\dots,0,\mp\cos\alpha)\in S^{d-1}$ be the unit vector corresponding to the angle  $\alpha\in[0,\frac{\pi}{2})$.
Then  
the center of mass of the  \textnormal{central section} $K\cap \xi^{\perp}$ is at the origin  if and only if
\begin{equation}\label{mihalna2}
({\mathcal R}(w_j\rho^{d}_{K}(w))(\xi)=0,\qquad j=1,\dots, d-1.
\end{equation}
Also, 
the  moments of inertia of  the central section  $K\cap \xi^{\perp}$ with respect to any $(d-2)$-dimensional subspace  $\Pi$ are  constant independent of $\Pi$ if and only if 
\begin{equation}\label{rabota1}
({\mathcal R}(w_1^2\rho^{d+1}_{K}(w))(\xi)=const \,(d+1)(1-\xi_1^2),
\end{equation}
\begin{equation}\label{mihalna1}
({\mathcal R}(w_j^2\rho^{d+1}_{K}(w))(\xi)=const\,(d+1)\qquad for \,\,all \quad j=2,\dots, d-1,
\end{equation}
and
\begin{equation}\label{mihalna111}
({\mathcal R}(w_jw_l\rho^{d+1}_{K}(w))(\xi)=0,\qquad j,l=1,\dots, d-1,\quad j\neq l.
\end{equation}
\eprop
\bp
If the center of mass of $K\cap\xi^{\perp}$ is at the origin, we have 
$$
\frac{1}{\textnormal{vol}_{d-1}(K\cap\xi^{\perp})}\int\limits_{K\cap\xi^{\perp}}xdx=0.
$$
Passing to  polar coordinates in $\xi^{\perp}$ 
and taking into account the fact that for $w\in \xi^{\perp}$ we have $w_d=w_1\tan\alpha $, 
we obtain the first statement of the lemma.

Let $\Pi$ be any $(d-2)$-dimensional subspace of $\xi^{\perp}$ and let $u=u_{d-1}$ be a unit vector in $\xi^{\perp}$ orthogonal to $\Pi$.  By (\ref{moment})  the condition on the moments reads as
\begin{equation}\label{dhob1}
I_{K\cap\xi^{\perp}}(\Pi)=\int\limits_{K\cap\xi^{\perp}}(x\cdot u)^2dx=const\qquad\quad\forall u\in S^{d-1}\cap\xi^{\perp}.
\end{equation}

Denote by $\iota_1,\dots\iota_{d-1}$  the orthonormal basis in $\xi^{\perp}$
such that $\iota_1=\cos\alpha\,\, e_1+\sin\alpha\, \,e_d$ and 
$\iota_j=e_j$ for $j=2,\dots, d-1$. 
Passing to polar coordinates and decomposing $u$ in the basis $\{\iota_j\}_{j=1}^{d-1}$, 
 we see that the  moments of inertia of  the central section  $K\cap \xi^{\perp}$ with respect to any $(d-2)$-dimensional subspace  are  constant if and only if 
 \begin{equation}\label{rabota2}
 ({\mathcal R}((w\cdot\iota_1)^2\rho^{d+1}_{K}(w))(\xi)=const\,(d+1),
\end{equation}
(\ref{mihalna1}) holds, 
and
 \begin{equation}\label{rabota3}
({\mathcal R}((w\cdot \iota_j)(w\cdot\iota_l)\rho^{d+1}_{K}(w))(\xi)=0,\qquad j,l=1,\dots, d-1,\quad j\neq l,
 \end{equation}
  (see the proof of 
Theorem 1 in \cite{R}). Since $w\cdot\iota_1=w_1\cos\alpha+w_d\sin\alpha$ and $w_d=w_1\tan\alpha$, we see that (\ref{rabota2}) and (\ref{rabota3}) are equivalent to (\ref{rabota1}) and (\ref{mihalna111}). This gives the second statement and the lemma is proved.
\ep

We remark  that for any body of revolution  around the $x_1$-axis, \eqref{mihalna2} holds for $j=2,\dots,d-1$.
Taking $u=\iota_j$  in the integral in (\ref{dhob1}), by rotation invariance we obtain   that
the moments in \eqref{mihalna1} are equal  for  $j=2,\dots,d-1$.
Also, arguing as at the end of  the proof of Lemma \ref{work1h} we  see  that
\eqref{mihalna111} is valid.

By these remarks, Step 2, Lemma \ref{cl1} with $s_o=1-3\tau$ and Proposition \ref{workHARD1} with $K=K_f$, when  $K_f$ is the  body of revolution we are constructing,
equations \eqref{mihalna2}, (\ref{rabota1}), \eqref{mihalna1} and  \eqref{mihalna111}    hold if $\tan\alpha\in (1-3\tau, 1-2\tau)$ with the constants in (\ref{rabota1}), \eqref{mihalna1} independent of $\xi$. 
Also,  the left-hand sides      of   (\ref{mihalna2}), (\ref{rabota1}) and (\ref{mihalna1})   are
already defined on the cap
$$
{\mathcal U}_{\tau}=\{\xi\in S^{d-1}:\,\xi_1=\pm\sin\alpha, \quad\alpha\in [0,\frac{\pi}{2}],\quad\tan\alpha \ge 1-3\tau   \}
$$
and are smooth even rotation invariant functions there.

Assume for a moment that we have constructed a smooth body $K_f$ so that conditions
\begin{equation}\label{mihalna33}
({\mathcal R}(w_1^2\rho^{d+1}_{K_f}(w))(\xi)=const \,(d+1)(1-\xi_1^2),\qquad 
({\mathcal R}(w_1\rho^{d}_{K_f}(w))(\xi)=0,
\end{equation}
hold for all unit vectors $\xi\in S^{d-1}$ with $\xi_1=\pm\sin\alpha$, corresponding to the angles  $\alpha\in[0,\frac{\pi}{2}]$ such that $\tan\alpha<1-2\tau$.
Then by the above remarks, Proposition \ref{workHARD1}  and the converse part of 
Lemma \ref{cl1} with $s_o=0$, 
conditions (\ref{ijn}), (\ref{cmassn}) of Proposition  \ref{hob1} are satisfied  for all $s>0$ and $K_f$ floats in equilibrium in every orientation at the level $\frac{\textrm{vol}_d(K)}{2}$.

Thus, it remains to construct the part of $K_f$ so that (\ref{mihalna33})
holds for all unit vectors $\xi$ corresponding to the angles  $\alpha\in[0,1-2\tau]$. 
To this end, denote by $\varphi_h$ and
$\psi_h$ the  left-hand sides of  (\ref{mihalna33}) defined on ${\mathcal U}_{\tau}$.
We put 
$\varphi_h(\xi)=$ $const\,(d+1)(1-\xi_1^2)$ and $\psi_h(\xi)=0$ 
for $\xi\in S^{d-1}$ such that $\xi_1=\pm\sin\alpha$ and
$\tan\alpha\in$ $ [0,1-2\tau]$. This definition
agrees with the one we already have when $\tan\alpha\in$ $ [1-3\tau,1-2\tau]$, so $\varphi_h$ and $\psi_h$ are even rotation invariant infinitely smooth functions
on the entire sphere.

\medskip

Recall that the values of ${\mathcal R}g(\xi)$ for all $\xi\in S^{d-1}$ such that $\xi_1=\pm\sin\alpha$ and $\tan\alpha$ $>1-3\tau$
are completely determined by the values of
the even  function $g(w)$ for all  $w\in S^{d-1}$ satisfying $w_1=\pm\cos\alpha$ and  $\tan\alpha>1-3\tau$. Moreover, for
 {\it bodies of revolution} (but not in general) the converse is also true (see the explicit inversion formula in  \cite[p. 433, formula (C.17)]{Ga}).

\medskip

Since  the equation  ${\mathcal R}g=\widetilde{g}$ with even $C^{\infty}$ right-hand side $\widetilde{g}$ is equivalent to
$$
\frac{g(\xi)+g(-\xi)}{2}={\mathcal R^{-1}}\widetilde{g}(\xi),
$$
we can rewrite  the equations in (\ref{mihalna33})     as
\begin{equation}\label{HOB1}
w_1^2(\rho^{d+1}_K(w)+\rho^{d+1}_K(-w))=2({\mathcal R^{-1}}\varphi_h)(w)
\end{equation}
and
\begin{equation}\label{HOB2}
w_1(\rho^{d}_K(w)-\rho^{d}_K(-w))=
2({\mathcal R^{-1}}\psi_h)(w).
\end{equation}
The already constructed part of $\rho_K$ satisfies these equations for the vectors $w\in S^{d-1}$ such that $w_1=\pm\cos\alpha$ and $\tan\alpha>1-3\tau$.

Since the spherical Radon transform commutes with rotations and our initial $\rho_K$ was rotation invariant,  the  even functions
$2{\mathcal R^{-1}}\varphi_h(w)$, $2{\mathcal R^{-1}}\psi_h(w)$
are rotation invariant as well  and can be written as $\Phi_h(\alpha)$ and $\Psi_h(\alpha)$, where $w\in S^{d-1}$ is such that $w_1=\pm\cos\alpha$ and $\alpha\in[0,\frac{\pi}{2}]$.
Note that the mappings $h\mapsto \Phi_h$, $h\mapsto\Psi_h$ are continuous from $C^{k+d}$ to $C^k$, say. Thus,  for all $h$ sufficiently close to zero in $C^{k+d}$,
$\Phi_h$ and $\Psi_h$ will be close to $\Phi_0\equiv2w_1^2$ and
$\Psi_0\equiv 0$  in $C^k$.

We will be looking for a rotation invariant solution $\rho_K$ of (\ref{HOB1}) and (\ref{HOB2}),  which will be described in terms of the two functions $R(\alpha)$ and $r(\alpha)$ related to it by (\ref{tryuk}).
Equations     (\ref{HOB1}) and (\ref{HOB2})    translate into
\begin{equation}\label{+r}
R^{d+1}(\alpha)+r^{d+1}(\alpha)=\frac{\Phi_h(\alpha)}{\cos^2\alpha},\quad R^{d}(\alpha)-r^{d}(\alpha)=\frac{\Psi_h(\alpha)}{\cos\alpha}.
\end{equation}
Equations (\ref{+r}), together with the conditions $R(\alpha)>0$ and $r(\alpha)>0$, determine $R(\alpha)$ and $r(\alpha)$ uniquely, and they coincide with the functions $R_h$ and $r_h$ obtained in Step 2
for  all $\alpha\in[0,\frac{\pi}{2}]$    with $\tan\alpha\ge 1-3\tau$. Therefore, any solution $R$,  $r$  of this system will satisfy $R(\alpha)=R_h(\alpha)$,
$r(\alpha)=r_h(\alpha)$ in this range.

If  $h$ and several of its derivatives are small enough, the functions $\Phi_h-2w_1^2$, $\Psi_h$ and several of their derivatives
are uniformly close to zero. 
Since  the map
$(R,r)\,\,\mapsto\,\,(
R^{d+1}+r^{d+1},\,
R^{d}-r^{d})$
is smoothly invertible near the point $(1,1)$  by the inverse function theorem,  the functions $R$, $r$ exist in this case on the entire interval $[0,\frac{\pi}{2}]$, and are close to $1$ in $C^2$. Moreover, $R'(0)=r'(0)=0$, because $\Phi^{'}_h(0)=0$,  $\Psi^{'}_h(0)=0$, (otherwise the functions ${\mathcal R^{-1}}\varphi_h$, ${\mathcal R^{-1}}\psi_h$ would not be smooth at $(1,0,\dots,0)$). This is enough to ensure that the body given by $R$ and $r$ is convex and corresponds to some strictly
concave function $f$ defined on $[-r(0),R(0)]$.

This completes the proof of Theorem \ref{bitsya} in the case of   odd dimensions. $\,\square$

\medskip

It remains to prove Theorem \ref{hob3}. Assume that
 a body $K\subset{\mathbb R^3}$ has   density ${\mathcal D}$ and  volume $V$. If $K$ 
is submerged in  liquid of density ${\mathcal D}'$ and 
$V'$ is the volume of a submerged part, then, by Archimedes' law, 
${\mathcal D}V={\mathcal D}'V'$, (cf. \cite[p. 257]{H}, \cite[p. 657]{Zh}). Taking ${\mathcal D}'=1$ and   $V'\!=\!\frac{1}{2}V$, 
we obtain the result. $\square$

\section{Appendix A:  proof of Theorem 
	\ref{olovzhal} from \cite{O}}

\subsection{The ``if$\,$" part}

We begin with several auxiliary lemmas. 

\bl\label{neozh1}
Let $d\ge 2$, let  $M\subset {\mathbb R^d}$ be a convex body and let $\varepsilon\in (0,1)$.  Consider the neighborhood of $\partial M$,
$U_{\varepsilon}$ $=U_{\varepsilon}(\partial M)=\{p\in {\mathbb R^d}:\, \textnormal{dist}(p,\partial M) <\varepsilon   \}$ and let $S(M)=S_{d-1}(M)$ be the $(d-1)$-dimensional surface area of $M$. Then
$\textnormal{vol}_d(U_{\varepsilon})\le 6\varepsilon S(M)$,
provided $\varepsilon$ is small enough.
\el
\bp
We fix a small $\varepsilon>0$ (we will choose it precisely later) and claim first that
\begin{equation}\label{neozh1}
\textrm{vol}_{d}(M\cap U_{\varepsilon} )\le \textrm{vol}_{d}(({\mathbb R^d}\setminus M)\cap U_{\varepsilon} ).
\end{equation}
Assume for a moment that $M$ is a convex polytope and consider the rectangular prisms $T_F$ based on facets $F$ of $M$ of height $2\varepsilon$, $T_F=F\times [-\varepsilon,\varepsilon]$ and such that $F\times (0,\varepsilon]\subset {\mathbb R^d}\setminus M$, $F\times [-\varepsilon,0]\subset  M$.
The  union of these prisms inside $M$
contains $M\cap U_{\varepsilon}$ and the parts of  prisms corresponding to the neighboring facets intersect. 
On the other hand,   the parts outside of $M$ do not intersect and the inequality for polytopes follows from
$$
\textrm{vol}_{d}(M\cap U_{\varepsilon} )\le\textrm{vol}_{d}\Big(\bigcup\limits_F (F\times [-\varepsilon,0])\Big)\le\qquad\qquad\qquad\qquad
$$
$$ \qquad\qquad\qquad\qquad\qquad\le\textrm{vol}_{d}\Big(\bigcup\limits_F (F\times [0,\varepsilon])\Big)\le \textrm{vol}_{d}(({\mathbb R^d}\setminus M)\cap U_{\varepsilon} ).
$$
 The general case can be obtained  by approximation of $M$ by polytopes and passing to the limit in the previous inequality. This proves the claim.

By (\ref{neozh1}) we have
$
\textnormal{vol}_d(U_{\varepsilon})\le 2 \textrm{vol}_{d}(({\mathbb R^d}\setminus M)\cap U_{\varepsilon} )
$
and it is enough to estimate the last volume. To do this, we will use the Steiner formula, \cite[p. 208]{Sch2}:
$$
\textnormal{vol}_d(M+\varepsilon B^d_2)=\sum\limits_{i=1}^{d}\varepsilon^{d-i}\kappa_{d-i}v_i(M),
$$
where 
$$
M+\varepsilon B^d_2=\{p=p_1+p_2\in{\mathbb R^d}:\,  p_1\in M\quad\textrm{and}\quad p_2\in \varepsilon B^d_2\},
$$
 and $v_i(M)$ are the intrinsic volumes of $M$, $1\le i\le d$, \cite[p. 214]{Sch2}. In particular, 
$v_d(M)=\textrm{vol}_d(M)$ and $v_{d-1}(M)$ is the surface area $S(M)$.
Since 
$$
({\mathbb R^d}\setminus M)\cap U_{\varepsilon} \,\,\subseteq \,\,(M+\varepsilon B^d_2)\setminus M,
$$
we obtain for $d=2$,
$$
\textrm{vol}_{2}(({\mathbb R^2}\setminus M)\cap U_{\varepsilon} )\le \sum\limits_{i=1}^{2}\varepsilon^{2-i}\kappa_{2-i}v_i(M)-\textrm{vol}_{2}(M)=\varepsilon \kappa_1v_1(M)=2\varepsilon S(M),
$$
 and for $d\ge 3$, 
$$
\textrm{vol}_{d}(({\mathbb R^d}\setminus M)\cap U_{\varepsilon} )\le \sum\limits_{i=1}^{d}\varepsilon^{d-i}\kappa_{d-i}v_i(M)-\textrm{vol}_{d}(M)=
$$
$$
2\varepsilon S(M)+\sum\limits_{i=1}^{d-2}\varepsilon^{d-i}\kappa_{d-i}v_i(M)\le 3\varepsilon S(M),
$$
provided $\varepsilon$ is so small that $\varepsilon (d-2)\max\limits_{1\le i\le d-2}(\kappa_{d-i}v_i(M))< S(M)$.
This gives the desired estimate.
\ep
To prove the next result we introduce some notation.
Let $P_H$ be the orthogonal projection onto a hyperplane $H$. For a small $\varepsilon>0$ we
let
$$
\Xi_{\varepsilon}=P_H(\{p\in \partial K:\,\textrm{dist}(p,H)<\varepsilon\}).
$$
Let $D$ be the length of a diameter of $K$ and let 
$\mu=\frac{2D^{d}}{\textrm{vol}_d(K\cap H^-(\xi))}$. We put
\begin{equation}\label{sig1}
 \Sigma_{\mu\varepsilon}=\{p\in H(\xi):\, \textrm{dist}(p, \partial K\cap H(\xi))<\mu\varepsilon\},
\end{equation}
where $H(\xi)$ is a hyperplane for which  (\ref{fubu}) holds.

\begin{figure}[ht]
	\includegraphics[width=320pt]{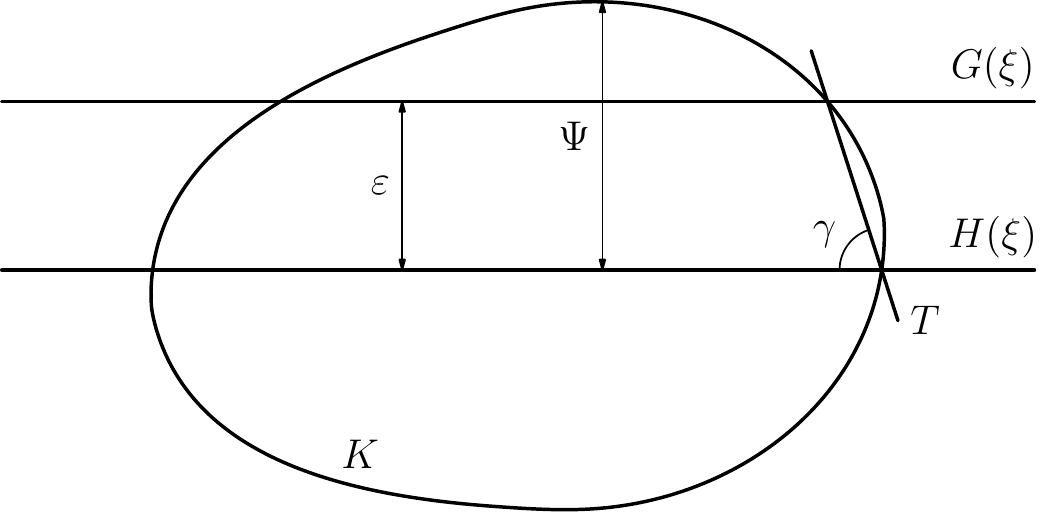}
	\caption{The hyperplanes  $H(\xi)$, $G(\xi)$, and  $T$.}
	\label{figO}
\end{figure}

\bl\label{sigma}
We have $\Xi_{\varepsilon}\subset \Sigma_{\mu\varepsilon}$, and 
$\textnormal{vol}_{d-1}(\Sigma_{\mu\varepsilon})<6c_d\mu D^{d-2}\varepsilon\to 0$ as $\varepsilon\to 0$.
\el
\bp
Consider a hyperplane $G(\xi)\in H^-(\xi)$ which is parallel to $H(\xi)$ and such that $\textrm{dist}(H(\xi), G(\xi))=\varepsilon$ for  $\varepsilon>0$ small enough. 
Consider also  a hyperplane $T$ containing any two corresponding parallel $(d-2)$-dimensional planes that support  $K\cap H(\xi)$ and $K\cap G(\xi)$. In the half-space $H^-(\xi)$ containing these sections choose an angle $\gamma$ between $T$ and $H(\xi)$ which is not obtuse (see Figure \ref{figO}, cf. Figure 1 in \cite{O}).

Denote by $\Psi$ the maximal distance between $H(\xi)$ and  any point in $K\cap H^-(\xi)$. 
Then
$$
\Psi\le D\sin\gamma,\qquad \textrm{vol}_d(K\cap H^-(\xi))<D^{d-1}\Psi\le D^{d}\sin\gamma.
$$
On the other hand, if 
$\lambda=\frac{\textrm{vol}_d(K\cap H^-(\xi))}{\textrm{vol}_d(K)}$, then
$$
\textrm{vol}_d(K\cap H^-(\xi))\ge\frac{\lambda}{1+\lambda}\textrm{vol}_d(K)\ge \frac{1}{2}\lambda\,\textrm{vol}_d(K),
$$
which yields 
$$
\sin\gamma>\frac{\lambda\textrm{vol}_d(K)}{2D^{d}},\qquad |\cot\gamma|<\frac{2D^{d}}{\lambda\textrm{vol}_d(K)}=\mu.
$$
Since the distance between the corresponding $(d-2)$-dimensional support planes to $K\cap H(\xi)$ and $P_{H(\xi)}(K\cap G(\xi))$ is $\varepsilon|\cot\gamma|<\mu\varepsilon$, 
we see that $\Xi_{\varepsilon}$ is a subset of  $\Sigma_{\mu\varepsilon}$.

Let $S$ be the $(d-2)$-dimensional surface area of $\partial K\cap H(\xi)$. Then
$$
\textrm{vol}_{d-1}(\Sigma_{\mu\varepsilon})\le 6\mu\varepsilon S(K\cap H(\xi))<6\mu \varepsilon c_dD^{d-2}\to 0,\quad\textrm{as}\quad \varepsilon\to 0.
$$
The first inequality follows from Lemma \ref{neozh1}, provided we identify $H(\xi)$ with ${\mathbb R}^{d-1}$  and put  $M=K\cap H(\xi)$. In the second inequality we used the fact that the surface area of $\partial K\cap H(\xi)$ does not exceed $c_dD^{d-2}$, where $c_d$ is some constant depending on the dimension, (it follows, for example,  from inequality (7) in  \cite[Theorem 1]{CSG}).
\ep

Now consider a family ${\mathcal W}={\mathcal W}_{\Gamma}$ of hyperplanes $H$ satisfying 
(\ref{fubu}) which are parallel to some $(d-2)$-dimensional subspace $\Gamma$. Each such hyperplane is determined by the angle $\theta\in[0,2\pi]$ it makes with some fixed $H_0\in {\mathcal W}$ (we take the orientation  into account). We will denote by $H(\theta)$ and $H(\theta+\varDelta\theta)$ the hyperplanes in ${\mathcal W}$ making angles $\theta$ and $\theta+\varDelta\theta$ with the chosen $H_0=H(0)=H(2\pi)$.
\bl\label{au2}
For sufficiently small $\varDelta\theta$ the $(d-2)$-dimensional plane 
$H(\theta)\cap H(\theta+\varDelta\theta)$ passes through $K$.
\el
\bp
Observe  first that for $\varDelta\theta$ small enough, the compact convex sets $K\cap H^-(\theta)$ and $K\cap H^-(\theta+\varDelta\theta)$ have a common point in the interior of $H^-(\theta)$.
Indeed, let $\beta$ be the smallest angle between $H(\theta)$ and the supporting hyperplanes to $K$ at points in $\partial K\cap H(\theta)$. As in the proof of Lemma \ref{sigma},
one can show that
$$
\beta>\sin\beta>\frac{\lambda \textrm{vol}_d(K)}{2D^{d}}=\frac{1}{\mu}.
$$
Therefore, any supporting hyperplane to $K$ making a positive angle  with $H(\theta)$ which is less than $\frac{1}{\mu}$,
must also support $K\cap H^-(\theta)$. 
Let  $\widetilde{H}(\theta+\varDelta\theta)$ be the supporting hyperplane 
to $K\cap H^-(\theta+\varDelta\theta)$
 parallel to $H(\theta+\varDelta\theta)$. Then 
$\widetilde{H}(\theta+\varDelta\theta)$ is also  the supporting hyperplane 
to $K\cap H^-(\theta)$, provided
$\varDelta\theta<\frac{1}{\mu}$. This proves  the observation.

Using the observation, we see that if  $H(\theta)\cap H(\theta+\varDelta\theta)$ does not pass through $K$, then
$K\cap H^-(\theta)$ and $K\cap H^-(\theta+\varDelta\theta)$ 
are contained in one another. This   contradicts the fact that
they 
have the same volume and the result follows.
\ep

Now choose a ``moving" system of coordinates in which the $(d-2)$-dimensional plane
$H(\theta)\cap H(\theta+\varDelta\theta)$ is the $p_1p_2\cdots p_{d-2}$-coordinate plane, the axis $p_{d-1}$ is in $H(\theta)$ and the axis $p_d$ is orthogonal to $H(\theta)$. We can assume that $\varDelta\theta$ is acute and is less than $\frac{1}{\mu}$.

The next lemma is a direct consequence of the fact that all hyperplanes in ${\mathcal W}$ satisfy (\ref{fubu}). Denote by $A\triangle B$ the symmetric difference of two sets $A$ and $B$, i.e., 
$A\triangle B=(A\setminus B)\cup (B\setminus A)$.

\begin{figure}[ht]
	\includegraphics[width=360pt]{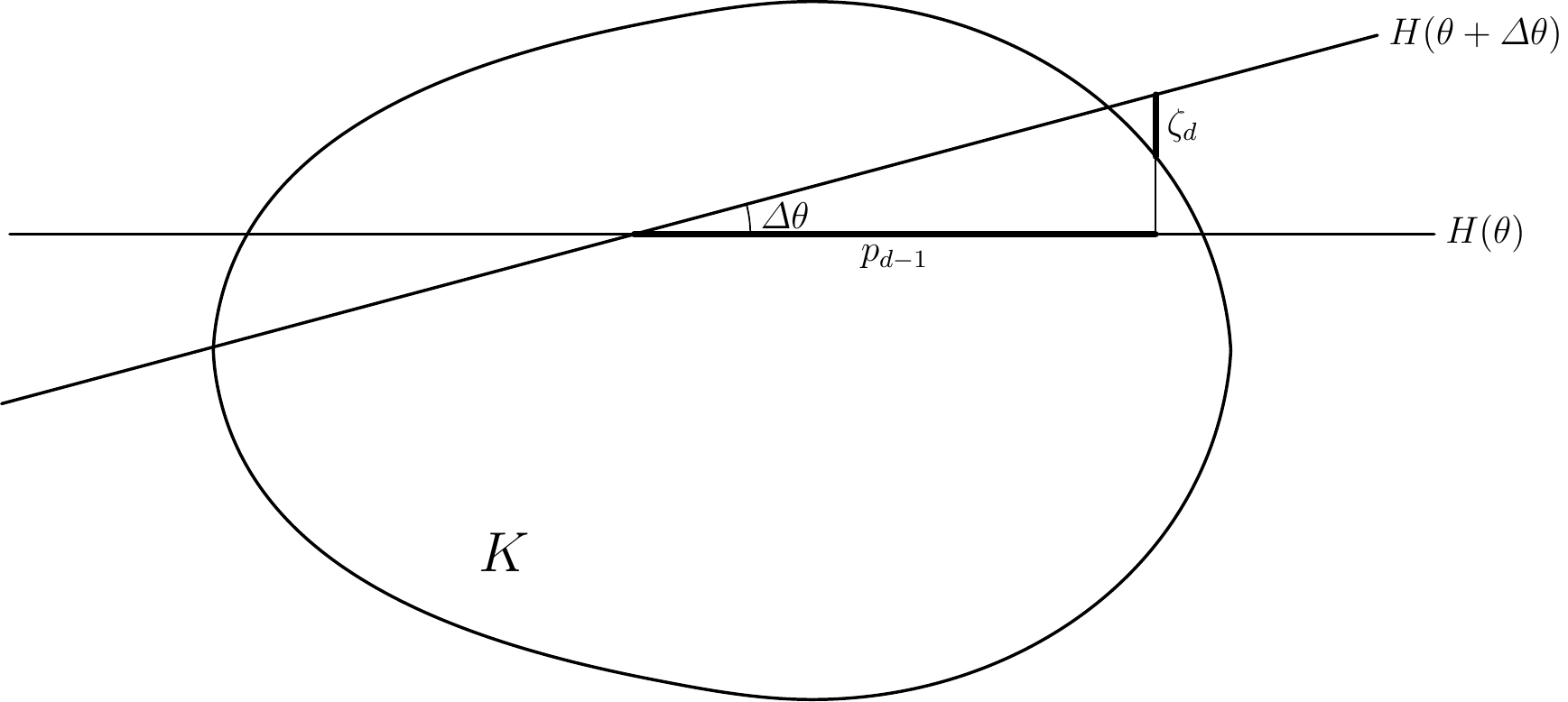}
	\caption{The  function $\zeta_d$.}
	\label{Fig7}
\end{figure}

\bl\label{tou1}
Let
$\Lambda= (K\cap H(\theta))\triangle P_{H(\theta)} (K\cap H(\theta+\!\varDelta\theta))$.
Then
\begin{equation}\label{kr1}
\varDelta V=\textnormal{vol}_d(K\cap H^-(\theta))-\textnormal{vol}_d(K\cap H^-(\theta+
\varDelta\theta))=
\end{equation}
$$
\int\limits_{K\cap H(\theta)}p_{d-1}\tan \varDelta \theta \,dp-
\int\limits_{\Lambda}\zeta_d \,dp=0,
$$
where  $p_{d-1}=p_{d-1}(\theta,\!\varDelta\theta)$ and 
$\zeta_d=\zeta_d(\theta,\!\varDelta\theta)$ is an error  of $p_d=p_{d-1}\tan\varDelta\theta$ in $\Lambda$ which is obtained  during the computation of $\varDelta V$ using the first integral above $\textnormal{(}$see Figure \ref{Fig7}$\textnormal{)}$. 
\el

We are ready to finish the  proof of the $``\,\textrm{if}\,"$ part  of Theorem \ref{olovzhal}. Let $p_{d-1}({\mathcal C}(K\cap H(\theta)))$ be the $(d-1)$-coordinate of ${\mathcal C}(K\cap H(\theta))$ with respect to the moving coordinate system.
By (\ref{kr1}), we have
$$
p_{d-1}({\mathcal C}(K\cap H(\theta)))=\frac{\int\limits_{K\cap H(\theta)}p_{d-1} \,\,\,dp}{\textrm{vol}_{d-1}(K\cap H(\theta))}=\frac{\int\limits_{\Lambda}\zeta_d \,\,\,dp}{\textrm{vol}_{d-1}(K\cap H(\theta))\tan \varDelta \theta}.
$$
Since for every $ p\in \Lambda$ there exists $ q\in \Xi_{D\sin\varDelta\theta}$  such that $P_{H(\theta)}q=p$,  applying Lemma \ref{sigma}
we see that 
$$
\textrm{vol}_{d-1}(\Lambda)\le  \textrm{vol}_{d-1}(\Xi_{D\sin\varDelta\theta})\le \textrm{vol}_{d-1}(\Sigma_{\mu D\sin\varDelta\theta})\le2c_d\mu D^{d-1}\varDelta\theta\to 0
$$
as $\varDelta\theta\to 0$.
Using the estimate $|\zeta_d|\le D\tan\varDelta\theta$, the previous inequalities and  the fact that $\Lambda\subset \Sigma_{\mu D\sin\varDelta\theta}$,
we obtain
$$
|p_{d-1}({\mathcal C}(K\cap H(\theta)))|\le \frac{D\tan \varDelta \theta\,\,\textrm{vol}_{d-1}(\Lambda)}{\textrm{vol}_{d-1}(K\cap H(\theta))\tan \varDelta \theta}\to 0
$$
as $\varDelta\theta\to 0$.
We see that, as $\varDelta\theta\to 0$,  the $(d-2)$-dimensional plane 
$H(\theta)\cap H(\theta+\varDelta\theta)$ tends to a limiting position   that passes through the center of mass of $K\cap H(\theta)$. 

To show that ${\mathcal C}(K\cap H(\theta))$ is the characteristic point of $H(\theta)$,
it is enough to take any $(d-2)$-dimensional subspace $\Gamma'$ that is parallel to $H(\theta)$, and to repeat the above considerations for the family  of hyperplanes ${\mathcal W}_{\Gamma'}$ that are parallel to $\Gamma'$.

Since the  subspace $\Gamma$ and the angle $\theta$ were chosen arbitrarily, we obtain the proof of the ``$\,\textrm{if}\,$" part of the theorem.

\medskip

\subsection{Proof of the converse part of Theorem \ref{olovzhal}}

Let $\Gamma$ be an arbitrary $(d-2)$-dimensional subspace and let ${\mathcal V}$ be a family of hyperplanes $H$ parallel to $\Gamma$ and such that  for all $H\in{\mathcal V}$ the centers of mass of $K\cap H$ coincide with the characteristic points of $H$.
Also, as above, choose an arbitrary angle $\theta$, the hyperplanes $H(\theta)$ and $H(\theta+\varDelta\theta)$ in ${\mathcal V}$ and  a ``moving" coordinate system. Since 
${\mathcal C}(K\cap H(\theta))$ is the characteristic point of $H(\theta)$
we can assume that  $p_{d-1}({\mathcal C}(K\cap H(\theta)))\to 0$ as $\varDelta\theta\to 0$. 

Using (\ref{kr1}) we have
$$
\frac{\varDelta V}{\varDelta \theta}=\frac{\tan \varDelta \theta}{\varDelta\theta}\,
\int\limits_{K\cap H(\theta)}p_{d-1}\,\,dp-
\int\limits_{\Lambda}\frac{\zeta_d}{\varDelta\theta} \,\,\,dp.
$$
Since 
${\mathcal C}(K\cap H(\theta+\varDelta\theta))\to {\mathcal C}(K\cap H(\theta))$ and $ \partial K\cap H(\theta+\varDelta\theta)\to$ $ \partial  K\cap H(\theta)$ as $\varDelta\theta\to 0$,
the set $\Lambda$ defined in Lemma \ref{tou1} satisfies $
\textrm{vol}_{d-1}(\Lambda)\to 0$ as $\varDelta\theta\to 0$. Using this and the fact that 
$|\zeta_d|\le D\tan\varDelta\theta$ we see that both summands in the right-hand side of the above identity tend to $0$ as $\varDelta\theta\to 0$. This gives $
\lim\limits_{\varDelta\theta\to 0}\frac{\varDelta V}{\varDelta\theta}=0$. 

 Now consider the function $\xi\mapsto g(\xi):=\textnormal{vol}_d(K\cap H^-(\xi))$ on $S^{d-1}$, where 
 $H(\xi)$ is the hyperplane from our family ${\mathcal V}$. By condition of the theorem, for every  $\xi\in S^{d-1}$ the center of mass ${\mathcal C}(K\cap H(\xi))$ is the characteristic point of $H(\xi)$ and     for any sequence $\{\xi_k\}_{k=1}^{\infty}$, $\xi_k\in S^{d-1}$, converging to  $\xi$  as $k\to\infty$
we have ${\mathcal C}(K\cap H(\xi_k))\to {\mathcal C}(K\cap H(\xi))$.
Since $\Gamma$ and $\theta$ were chosen arbitrarily, 
writing $g(\xi)$ in terms  of the  spherical angles $\varphi_1,\dots, \varphi_{d-1}$, $\varphi_j\in[0,\pi)$, $j=1$, $\dots,d-2$, $\varphi_{d-1}\in [0,2\pi)$, we can choose the corresponding sequences 
$\{\xi_{j, k}\}_{k=1}^{\infty}$, $\xi_{j, k}\in S^{d-1}$, converging to $\xi(\varphi_1,\dots, \varphi_{d-1})$, so that 
$\frac{\partial }{\partial \varphi_j}g(\xi)=0$ for all $\xi\in S^{d-1}$ and all $j=1,\dots, d-1$.
 Therefore,  $g$  must  be constant on $S^{d-1}$.
The proof of the converse  part is complete.

This finishes the proof of Theorem \ref{olovzhal}.$\qquad\qquad\qquad\qquad\qquad\qquad\qquad\quad\square$

\section{Appendix B:  proof of the converse part of Theorem 
	\ref{Fedja1}}

We start by recalling the so-called First Theorem of Dupin, 
(cf. \cite[pp. 658-660]{Zh}  and \cite[pp. 275-279]{DVP}; see also  \cite[Theorem 4]{R}). 

It was proved in \cite[Theorem 1.2]{HSW} that the surface of centers ${\mathcal S}$ is  $C^{k+1}$-smooth, provided $K$ is of class $C^k$, $k\ge 0$.
In particular,  if $K$ is an arbitrary convex body then  ${\mathcal S}$ is  $C^1$-smooth.

Let $\Gamma$  be any  $(d-2)$-dimensional subspace of ${\mathbb R^d}$. We   let  the family ${\mathcal W}={\mathcal W}_{\Gamma}$ of hyperplanes $H(\theta)$, $\theta\in[0,2\pi]$,  satisfying 
(\ref{fubu}) and which are parallel to $\Gamma$  be as in the previous section. We will use the notation  ${\mathcal C}(\theta)\in {\mathcal S}$ for the centers of mass of the corresponding ``submerged" parts $K\cap H^-(\theta)$ and
${\mathcal H}(\theta)$ for the tangent hyperplane to ${\mathcal S}$ at ${\mathcal C}(\theta)$.

\bt\label{D1}
Let $d\ge 2$, let $K\subset {\mathbb R^d}$ be a convex body  and let $\delta\in (0,\textnormal{vol}_d(K))$.  Then for any $\Gamma$ and for any $H(\theta)\in {\mathcal W}_{\Gamma}$, $\theta\in [0,2\pi]$, 
 ${\mathcal H}(\theta)$
is parallel to 
$H(\theta)$. Also, the bounded set
$L=L({\mathcal S})$ with  boundary 
${\mathcal S}$ is a strictly convex body.
\et
\bp
Fix $\Gamma$ and $\theta\in [0,2\pi)$.
 Rotating and translating if necessary  we can assume that   $H(\theta)=e_d^{\perp}$ and  $K\cap H^-(\theta)\subset\{p\in {\mathbb R^d}:\,p_d\le 0\}$.  
 Let $H(\widetilde{\theta})\in  {\mathcal W}_{\Gamma}$, $\widetilde{\theta}\neq\theta$, $\widetilde{\theta}\in[0,2\pi)$.  
 We claim  that ${\mathcal C}(\widetilde{\theta})$ is ``above" ${\mathcal C}(\theta)$, i.e., $p_d({\mathcal C}(\theta))<p_d({\mathcal C}(\widetilde{\theta}))$. 
Indeed, since  $p_d>0$
$\forall p\in (K\cap H^-(\widetilde{\theta}))\setminus  (K\cap H^-(\theta))$ but $p_d\le 0$ 
$\forall p\in (K\cap H^-(\theta))\setminus  (K\cap H^-(\widetilde{\theta}))$, we have
$$
p_d({\mathcal C}(\theta))=\frac{1}{\delta}\Big( \int\limits_{(K\cap H^-(\theta))\setminus  (K\cap H^-(\widetilde{\theta}))}p_d dp+\int\limits_{K\cap H^-(\theta)\cap H^-(\widetilde{\theta})}p_d dp\Big)<
$$
$$
\frac{1}{\delta}\Big( \int\limits_{(K\cap H^-(\widetilde{\theta}))\setminus  (K\cap H^-(\theta))}p_d dp+\int\limits_{K\cap H^-(\theta)\cap H^-(\widetilde{\theta})}p_d dp\Big)=\,p_d({\mathcal C}(\widetilde{\theta}))
$$
and  the claim is proved.

Now let ${\mathcal H}(\theta)$ be the hyperplane passing through ${\mathcal C}(\theta)$ which is parallel to $H(\theta)$ and let
${\mathcal H}^{\pm}(\theta)$ be the corresponding half-spaces.
Since $\Gamma$ and $\theta$ were chosen arbitrarily, we see that  ${\mathcal S}\subset {\mathcal H}^+(\theta)$. Since ${\mathcal S}$ is $C^1$-smooth and ${\mathcal S}\cap {\mathcal H}(\theta)={\mathcal C}(\theta)$, the hyperplane 
${\mathcal H}(\theta)$ is tangent to ${\mathcal S}$ at ${\mathcal C}(\theta)$.

Thus, for any $\xi\in S^{d-1}$  we have
${\mathcal S}\subset {\mathcal H}^+(\xi)$,
${\mathcal S}\cap {\mathcal H}(\xi)={\mathcal C}_{\delta}(\xi)$ and $\min\limits_{\{\xi\in S^{d-1}\}}|{\mathcal C}(K)-{\mathcal C}_{\delta}(\xi)|>0$. We conclude that
$L({\mathcal S})=\bigcap\limits_{\{\xi\in S^{d-1}\}}{\mathcal H}^+(\xi)$ is a strictly convex body.
\ep

To prove the converse part of Theorem 
\ref{Fedja1}
it is enough to show that   the orthogonal projection of ${\mathcal S}$ onto any $2$-dimensional subspace of ${\mathbb R^d}$ is a disc. Indeed, by
applying  \cite[Corollary 3.1.6, p. 101]{Ga} to $L({\mathcal S})$, we obtain that in this case ${\mathcal S}$ is a sphere. 
Using Theorem \ref{D1}, as well as the fact that all normal lines of the sphere intersect at its center,  we see that for  every $\xi\in S^{d-1}$ the lines $\ell(\xi)$ passing through ${\mathcal C}(K)={\mathcal C}({\mathcal S})$ and ${\mathcal C}_{\delta}(\xi)$ are orthogonal to $H(\xi)$. 
By  Definition \ref{efb} this means that $K$ floats in equilibrium in every orientation.

Let $\Gamma$ be as above, let $\Gamma^{\perp}$  be the $2$-dimensional subspace orthogonal to $\Gamma$ and let $P=P_{\Gamma^{\perp}}$ be the orthogonal projection  
onto $\Gamma^{\perp}$. 
To show that  $P({\mathcal S})$ is a disc for every $\Gamma$,  
we will prove the following lemma.
\bl\label{FO1}
 Let 
 $\xi(\theta)\in S^{d-1}$   be the normal vector to $H(\theta)$, let $\beta$ be a  closed curve  $\{{\mathcal C}(\theta):\,\theta\in[0,2\pi]\}\subset {\mathcal S}$ and let  $P\beta=\{P{\mathcal C}(\theta):\,\theta\in[0,2\pi]\}$ be  parametrized as $\theta\mapsto \varrho(\theta)$, $\theta\in [0,2\pi]$.
Then
\begin{equation}\label{FL1}
\varrho'(\theta)=-\frac{1}{\delta}I_{K\cap H(\theta)}(\Pi)\,\,\xi'(\theta)\qquad \forall \theta\in [0,2\pi],
\end{equation}
where $\Pi$ is the $(d-2)$-dimensional plane passing through ${\mathcal C}(K\cap H(\theta))$ and parallel to $\Gamma$.
\el
Assume for a moment that (\ref{FL1}) is proved. 
By conditions of the theorem,  $I_{K\cap H(\theta)}(\Pi)$ is the constant $c$  independent of $\Pi$ and $\theta$. Integrating  both parts in (\ref{FL1}) we have $\varrho(\theta)=-c \,\xi(\theta)+C$, where  $C$ is  a constant vector. Hence,   $P\beta$ is a circle. Since 
 $\Gamma$ was chosen arbitrarily, the projection of ${\mathcal S}$ onto any $2$-dimensional subspace is a disc.

To finish the proof, it remains to prove 
 the lemma. 
 \bp
 We can assume that  $H(\theta)=e_d^{\perp}$, $K\cap H^-(\theta)\subset \{p\in {\mathbb R^d}:\,p_d\le 0\}$ and  $\rho(\theta)$, $\xi(\theta)$, $\xi'(\theta)$ are $2$-dimensional, i.e.,  $\varrho(\theta)= (\varrho_{d-1}(\theta),\varrho_d(\theta))$, $\xi(\theta)=(0,1)$, $\xi'(\theta)=(-1,0)$.
Since the tangent vector $\varrho'(\theta)$ is parallel to $ {\mathcal H}(\theta)$ and since ${\mathcal H}(\theta)$  is parallel to $H(\theta)$ by the previous theorem, we conclude that $\varrho_d'(\theta)=0$. 

To compute $\varrho_{d-1}'(\theta)$, we will estimate  $\varrho_{d-1}(\theta+\varDelta\theta)-\varrho_{d-1}(\theta)$ for  $\varDelta\theta$ small enough. 
As in the previous appendix, we choose a ``moving" system of coordinates in which the $(d-2)$-dimensional plane
$H(\theta)\cap H(\theta+\varDelta\theta)$ is the $p_1p_2\cdots p_{d-2}$-coordinate plane. 
We  have
$$
\varrho_{d-1}(\theta+\varDelta\theta)-\varrho_{d-1}(\theta)=
\frac{1}{\delta}\Big( \int\limits_{K\cap H^-(\theta+\varDelta\theta)} p_{d-1}dp -\int\limits_{K\cap H^-(\theta)} p_{d-1}dp \Big)=
$$
$$
=\frac{1}{\delta}\Big( \int\limits_{K\cap H(\theta)}p^2_{d-1}\tan \varDelta \theta \,dp-
\int\limits_{\Lambda}p_{d-1}\zeta_d \,dp\Big),
$$
where the last equality is similar to (\ref{kr1}),
$\Lambda$ and $\zeta_d$ are as in Lemma \ref{tou1}  $\textnormal{(}$see Figure \ref{Fig7}$\textnormal{)}$.
Dividing both parts  by $\varDelta\theta$, passing to the limit as $\varDelta\theta\to 0$ and using the ``$\,$if$\,$" part of the theorem proved in the previous appendix, we obtain
$$
\varrho_{d-1}'(\theta)=\frac{1}{\delta}\int\limits_{ K\cap  H(\theta)}p^2_{d-1}dp=\frac{1}{\delta}I_{K\cap H(\theta)}(\Pi).
$$
This gives (\ref{FL1}).
\ep


{\bf Acknowledment}. The author is very thankful  to Mar\'ia de los \'Angeles Alfonseca Cubero, Mar\'ia de los \'Angeles Hern\'andez Cifre, Alexander Fish, Fedor Nazarov, Alina Stancu, Peter V\'arkonyi  and  Vlad Yaskin  for  their invaluable help and very useful discussions.  He  is also very indebted to the anonymous referees  whose important  remarks helped to improve the paper.

\end{document}